%% file: main.tex
\documentclass[onefignum,onetabnum]{siamart190516}
\input{ex_shared}

\ifpdf
\hypersetup{
  pdftitle={P2NN},
  pdfauthor={M.Liu, Z.Cai}
}
\fi

\DeclareUnicodeCharacter{2212}{-}
\usepackage[utf8]{inputenc}
\usepackage{bm}
\usepackage{mathtools}
\usepackage{indentfirst}
\usepackage{subfigure}
\usepackage{wrapfig}
\usepackage{tcolorbox}
\usepackage{color}
\usepackage[export]{adjustbox}
\usepackage{makecell}
\usepackage{textcomp}
\usepackage{algorithm }
\usepackage{algorithmic}
\usepackage{booktabs}
\usepackage{multirow}

\newcommand{\R}{\mathbb{R}}

\usepackage{graphicx}
\usepackage{diagbox}
\usepackage{pifont}
\newcommand{\vertiii}[1]{{\left\vert\kern-0.25ex\left\vert\kern-0.25ex\left\vert #1 
    \right\vert\kern-0.25ex\right\vert\kern-0.25ex\right\vert}}

\newcommand{\argmin}{\mathop{\mathrm{arg\,min}}}
    
\renewcommand{\div}{\nabla \cdot}

\newcommand{\bSigma}{\mbox{\boldmath${\Sigma}$}}

\newcommand{\bTheta}{\mbox{\boldmath${\Theta}$}}

\newcommand{\btau}{\mbox{\boldmath${\tau}$}}

\newcommand{\bvarphi}{\mbox{\boldmath${\varphi}$}}

\newcommand{\bomega}{\mbox{\boldmath${\omega}$}}
\newcommand{\bxi}{\mbox{\boldmath$\xi$}}

\newcommand{\bff}{{\bf f}}
\newcommand{\bF}{{\bf F}}
\newcommand{\bG}{{\bf G}}
\newcommand{\bH}{{\bf H}}
\newcommand{\bg}{{\bf g}}
\newcommand{\bn}{{\bf n}}

\newcommand{\bb}{{\bf b}}

\newcommand{\bx}{{\bf x}}
\newcommand{\bz}{{\bf z}}

\newcommand{\cI}{{\cal I}}

\newcommand{\cM}{{\cal M}}
\newcommand{\cQ}{{\cal Q}}
\newcommand{\cS}{{\cal S}}
\newcommand{\cT}{{\cal T}}

\newcommand{\cV}{{\cal V}}

\def\bp{{\bf p}}

\def\bb{{\bf b}}
\def\bc{{\bf c}}

\def\br{{\bf r}}
\def\bx{{\bf x}}
\def\by{{\bf y}}

\def\cM{{\cal M}}

\def\cP{{\cal P}}
\def\cS{{\cal S}}
\def\cT{{\cal T}}

\def\cP{{\cal P}}
\def\p2nn{{$\text{P}^2$NN}}

\def\div{{\mbox{\bf div\,}}}
\def\divt{${\mbox{\bf div}}_{_\cT}$\,}

\newcommand{\jump}[1]{[\![ #1]\!]}

\begin{document}

\maketitle
\begin{abstract}
This chapter offers a comprehensive introduction to the least-squares neural network (LSNN) method introduced in \cite{Cai2021linear, Cai2023nonlinear}, for solving scalar first-order hyperbolic partial differential equations, specifically linear advection-reaction equations and nonlinear hyperbolic conservation laws. The LSNN method is built on an equivalent least-squares formulation of the underlying problem on an admissible solution set that accommodates discontinuous solutions. It employs ReLU neural networks (in place of finite elements) as the approximating functions, uses a carefully designed physics-preserved numerical differentiation, and avoids penalization techniques such as artificial viscosity, entropy condition, and/or total variation. This approach captures shock features in the solution without oscillations or overshooting. Efficiently and reliably solving the resulting non-convex optimization problem posed by the LSNN method remains an open challenge. %\textcolor{red}{
This chapter concludes with a brief discussion on application of the structure-guided Gauss-Newton (SgGN) method developed recently in \cite{SgGN} for solving  shallow NN approximation.%}
\end{abstract}
%This chapter aims to introduce neural network-based discretization approaches, with a specific focus on the LSNN approach for hyperbolic conservation laws. The LSNN method stands out as an innovative technique within this field. To better position LSNN within the broader landscape of discretization methods, it would be beneficial to include an overview of various neural network-based techniques. This broader context will offer readers a clearer understanding of where LSNN fits among contemporary approaches.

%Moreover, to enhance the accessibility of the content for students, the chapter should allocate more space and attention to detailed explanations. By incorporating illustrative examples, step-by-step derivations, and intuitive analogies, the material can become more comprehensive and easier to grasp for learners at different levels of familiarity with the subject.

%The resulting discretization of the LSNN method is a non-convex optimization problem in the NN parameters. This high dimensional, non-convex optimization problem tends to be computationally intensive and complicated, and has been a bottleneck in using NNs to numerically solve PDEs.


\section{Introduction}\label{s:introduction}

Over the past five decades, numerous advanced mesh-based numerical methods have been developed for solving nonlinear hyperbolic conservation laws (HCLs) (see, e.g., \cite{leveque1992numerical, GoRa:96, shu1998essentially, Le:02, thomas2013numerical, hesthaven2017numerical}). However, accurately approximating solutions to HCLs remains computationally challenging due to two key difficulties. First, the location of the discontinuities in the solution is typically unknown in advance. Second, the strong form of the partial differential equation (PDE) becomes invalid at points where the solution is discontinuous. 

Recently, neural networks (NNs) have emerged as a novel class of functions approximators for solving partial differential equations (PDEs) (see, e.g., \cite{CaChLiLi:20,Weinan18,raissi2019physics,Sirignano18} and Section~\ref{NN_D}). A neural network function is a linear combination of compositions of linear transformations and a nonlinear univariate activation function. %\textcolor{red}{
One commonly used activation function is the Rectified Linear Unit (ReLU), defined as $\sigma(s) = \max\{0, s\}$. A ReLU neural network thus constructs a function as a piecewise linear approximation, with ``break points'' determined by its parameters, effectively forming a data-adaptive partition of the domain. %} 
As demonstrated in \cite{Cai2021linear, Cai2023linear, CCL2024}, ReLU NNs can approximate discontinuous functions with unknown interfaces far more effectively than traditional approximating functions, such as polynomials or continuous/discontinuous piecewise polynomials defined on a quasi-uniform, predetermined mesh. This makes ReLU NNs particularly suitable for addressing the first challenge.  

The strong form of a hyperbolic PDE is typically written with partial derivatives along coordinate directions, supplemented by the Rankine-Hugoniot (RH) jump condition at discontinuity interfaces for HCLs. 
Due to the unknown location of these interfaces, enforcing the RH condition in computations is difficult, if not impossible. To address this, we reformulate the PDE using physically meaningful derivatives, allowing the new form of the PDE to remain well-defined at the interface (see \cref{pde1b} for the directional derivative and \cref{pde2b} for the divergence operator). By applying the $L^2$ least-squares principle to this reformulated PDE, we derive an equivalent least-squares minimization problem on an admissible solution set that accommodates discontinuous solutions. Through appropriate numerical integration for the integral and physics-preserved numerical differentiation for the physically meaningful derivative, the least-squares neural network (LSNN) method is established as minimizing the discrete counterpart of the least-squares functional over the set of NN functions.  

%To use NNs (instead of finite elements) to approximate the solution of hyperbolic PDEs, we proposed the space-time least-squares neural network (LSNN) method in \cite{Cai2021linear} for linear advection-reaction equations and in \cite{Cai2023nonlinear} for scalar nonlinear hyperbolic conservation laws. To derive the LSNN method, we first rewrite the underlying PDE by introducing physics meaningful differential operator so that the PDE is also well-defined at where the solution is discontinuous (see \cref{pde1} for the directional derivative and \cref{pde2b} for the divergence operator). The underlying problem is then reformulated as an equivalent least-squares minimization problem over a proper solution space that permits discontinuous solution. Approximating integral by proper numerical integration and the differential operator by physics-preserved discrete differential operator, the LSNN method is finally defined as to minimize the discrete counterpart of the least-squares functional over the set of NN functions.    

Without relying on penalization techniques such as inflow boundary conditions, artificial viscosity, entropy conditions, or total variation constraints, the LSNN method introduced in \cite{Cai2021linear, Cai2023nonlinear} effectively captures the shock of the underlying problem without oscillations or overshooting. Additionally, the LSNN method is substantially more efficient in terms of degrees of freedom (DoF) compared to adaptive mesh refinement (AMR) methods, which locate the discontinuity interface through an adaptive mesh refinement process.

Despite the impressive approximation capabilities of NNs, the discretization resulting from NN-based methods leads to a non-convex optimization problem in the NN parameters. This high-dimensional, non-convex optimization is often computationally intensive and complex, presenting a significant bottleneck in using NNs for numerically solving PDEs. Nonetheless, considerable research efforts are underway, with some promising progress in developing efficient and reliable iterative solvers (training algorithms) and in designing effective initializations \cite{SgGN, CaiDokFalHer2024a, CaiDokFalHer2024b}.

The chapter is organized as follows. \Cref{s:problem} describes the advection-reaction equation and the scalar nonlinear HCL, their equivalent least-squares formulations, and preliminaries.
ReLU neural network and its approximation property to discontinuous functions are introduced in \Cref{s:NN}. The physics-preserved numerical differentiation and the LSNN method are defined in \Cref{s:LSNN}. \Cref{s:solver} discusses efficient iterative solvers. Finally, numerical results for various benchmark test problems are given in \Cref{s:NE}.

% \subsection{\textcolor{red}{NN-Based Discretization Methods  and Numerical Issues
% }\label{NN_D}}
\subsection{Methodological Remarks and Related Work}\label{NN_D}

Since NN functions are inherently nonlinear with respect to certain parameters, it is both convenient and natural to discretize a PDE by reformulating it as an optimization problem through either natural energy minimization or manufactured least-squares (LS) principles. Consequently, existing NN-based numerical methods for solving PDEs fall into two main categories: (1) Energy-based methods, such as deep Ritz and finite neuron methods \cite{Weinan18, Xu2020, LiuCai2, LiCaRa23}, which employ the Ritz formulation for the primary variable \cite{EkTe:76}, and the dual neural network (DuNN) \cite{LIU2025}, which uses complementary energy minimization for the dual variable \cite{BrFo:91, EkTe:76}. (2) Deep LS methods, which formulate PDE residuals into manufactured least squares using various norm choices and problem forms \cite{Dissanayake94, Berg18, raissi2019physics, Sirignano18, CaChLiLi:20, Cai2021linear}. Energy-based approaches are applicable to a class of self-adjoint and positive definite problems that commonly arise in continuum mechanics. These methods are not only physics informed but also physics preserved, when using NN as approximating functions, it is natural to discretize the underlying problem based on an energy formulation. 

However, many scientific and engineering problems are non-self-adjoint and thus lack a natural minimization principle. In such cases, one can then apply the LS principle to create a manufactured one.
The LS principle offers great flexibility: all consistent equations and data can be incorporated into the LS functional. Yet, balancing the various terms is challenging; ill-scaled LS functionals can lead to suboptimal accuracy and inefficient training. The minimal requirement for a viable LS formulation is its equivalence to the original PDE; 
otherwise, the physical fidelity of the model is compromised.
% For a given PDE, there are many least-squares (LS) formulations depending on the $L^2$-based norms, forms of the PDE, etc. (see, e.g., monographs \cite{bg5, Jia:98}, articles \cite{clmm, CaLeWa:04, CaSt:04}, and references therein). The minimum criterion for a LS formulation is its equivalence to the underlying problem. Obviously, if a LS formulation is not equivalent, it then violates the underlying physics and would not lead to a viable numerical method. The LS principle is so flexible that one may use all consistent equations and/or data to form a LS functional. Nevertheless, it is very difficult to keep all terms of a LS functional in the same scale. Without a well-balanced LS functional, the resulting numerical method is in general not optimal in accuracy and encounters difficulties in training. 

A classical example is the Bramble–Schatz LS (BSLS) formulation \cite{bs1} for scalar elliptic PDEs, which applies the LS principle directly to the {\it strong} form of the PDE, the boundary and the initial conditions using appropriate Sobolev norms \cite{BrSc:94}. When the solution is sufficiently smooth, e.g., belongs to $H^2(\Omega)$ (the collection of functions whose second-order weak derivatives are square-integrable), the BSLS is equivalent to the underlying problem and yields a well-balanced formulation. However, it fails in the presence of singularities, such as those caused by geometric corners or material interfaces. To address this, LS methods have been extended to first-order systems with properly chosen norms, resulting in numerous viable LS finite element methods. These methods, grounded in $L^2$-based norms, have been thoroughly developed and analyzed for problems such as convection-diffusion-reaction, elasticity, and Stokes equations \cite{BoGu94, BoGu95, BoCMM98, BrLP97, BrLP01, CaMaMc:97, CaMaMc:97b}. In particular, the LS methods introduced in \cite{cai1994first, CaLeWa:04, CaSt:04} are based on the original physical first-order systems, ensuring equivalence with the underlying PDEs. When adapted for NN discretizations, these formulations preserve physical laws provided that the numerical differentiation of certain differential operators also respects the underlying physics.

For scalar nonlinear hyperbolic conservation laws, several NN-based numerical methods have been recently introduced by various researchers (\cite{PNAS2019, Cai2021nonlinear, Cai2021linear, Fuks20, raissi2019physics, Patel22}). Those methods can be categorized as the physics-informed neural networks (PINNs) \cite{PNAS2019, Fuks20, raissi2019physics, Patel22} and the LSNN methods \cite{Cai2021nonlinear, Cai2021linear, Cai2023nonlinear}. Both methods are based on the least-squares principle, but the former uses the discrete $l^2$ norm while the latter uses the continuous $L^2$ norm. This difference explains why the latter can employ ReLU activation function, that is not differentiable pointwisely, and adaptively find accurate numerical integration (see \Cref{s:LSNN}). Fundamentally, the former is based on the strong form (\cref{pde1} or \cref{pde2}) that {\it violates} the underlying physics (see \Cref{s:problem}), and the latter builds on a correct (``weak'') form (\cref{pde1b} or \cref{pde2b}) that {\it preserves} the underlying physics.

Due to such a violation, the original PINN produces unacceptable approximate solution of the underlying problem by scientific computing standard. %\textcolor{red}{by scientific computing standard}. 
This phenomenon was observed by several researchers, e.g., \cite{Fuks20, Patel22}. Furthermore, \cite{Fuks20} modified the loss function by penalizing the artificial viscosity term.
\cite{Patel22} applied the discrete $l^2$ norm to the boundary integral equations over control volumes instead of the differential equations over points and modified the loss function by penalizing the entropy, total variation, and/or artificial viscosity. Even though the least-squares principle permits freedom of various penalizations, choosing proper penalization constants can be challenging in practice and it affects the accuracy,  efficiency, and stability of the method. In contrast, the LSNN does not require any penalization constants (see \Cref{s:LSNN}).

%Second, the differential operator of the underlying problem is approximated by either automatic differentiation or standard finite difference quotient for the PINN and by specially designed discrete differential operator for the LSNN. For example, the LSNN uses discrete directional differential operator in \cite{Cai2021linear} for linear advection-reaction problems, and various traditional conservative schemes in \cite{Cai2021nonlinear} or discrete divergence operator in this paper (see \cite{cai21} for its first version) for nonlinear scalar hyperbolic conservation laws. }

%\textcolor{red}{The original PINN has limitations that have been addressed in several studies (see, e.g., \cite{Fuks20, Patel22}). For nonlinear scalar hyperbolic conservation laws,  \cite{Fuks20} found that the PINN fails to provide reasonable approximate solution of the PDE and modified the loss function by penalizing the artificial viscosity term. \cite{Patel22} applied the discrete $l^2$ norm to the boundary integral equations over control volumes instead of the differential equations over points and modified the loss function by penalizing the entropy, total variation, and/or artificial viscosity. Even though the least-squares principle permits freedom of various penalizations, choosing proper penalization constants can be challenging in practice and it affects the accuracy,  efficiency, and stability of the method. In contrast, the LSNN does not require any penalization constants. 

\section{Scalar Hyperbolic Partial Differential Equations}\label{s:problem}

Let $\Omega$ be a bounded open domain in ${\R}^d$ ($d=1, \,2$, or $3$) with Lipschitz boundary, and $I=(0, T)$ be the temporal interval. This section describes linear advection-reaction equations defined on $\Omega$ and scalar nonlinear hyperbolic conservation laws defined on $\Omega\times I$ and their equivalent least-squares formulations. 

\subsection{Advection-Reaction Equations}\label{s:arp}

Let $\bm{\beta}(\bx) = (\beta_1, \cdots, \beta_d)^t\in C^1(\bar{\Omega})^d$ be the advective velocity field having no stagnation point and $\gamma \in C(\bar{\Omega})$ be the reaction coefficient. Without loss of generality, assume that the magnitude of $\bm{\beta}(\bx)$ is one in $\Omega$, i.e., $|\bm{\beta}(\bx)|\equiv 1$.
Otherwise, the equation below in \cref{pde1} may be rescaled by dividing $|\bm{\beta}(\bx)|$.

Let $f \in L^2(\Omega)$ and $g \in L^2(\Gamma_-)$ be given scalar-valued functions, where $\Gamma_-$ is the inflow part of the boundary $\Gamma=\partial \Omega$ given by
\[
\Gamma_- = \{\bx\in\Gamma :\, \bm{\beta}(\bx)\! \cdot\! \bm{n}(\bx) <0\}
\]
with $\bm{n}(\bx)$ the unit outward normal vector to $\Gamma$ at $\bx\in \Gamma$.
Consider the following linear advection-reaction equation \cite{leveque1992numerical}
\begin{equation}\label{pde1}
    \left\{\begin{array}{rccl}
    %\bm{\beta}\cdot \!\nabla u 
    %u_{\bm\beta} 
    \sum\limits_{i=1}^d\beta_i(\bx) \dfrac{\partial u(\bx)}{\partial x_i} + \gamma u &= & f &\text{ in }\,\, \Omega, \\[2mm]
    u&=&g &\text{ on }\,\, \Gamma_{-}.
    \end{array}\right.
\end{equation}
If the inflow boundary data $g$ is discontinuous, so is the solution $u(\bx)$. Hence, the solution has no partial derivative along any coordinate direction on the discontinuity interface $\cI$, which, in turn, implies that the PDE in \cref{pde1} is not valid on the $\cI$. 

To deal with this issue, let us introduce characteristic curves of \cref{pde1} defined by
\begin{equation}\label{l-chara}
    % \dfrac{d\bx_{\subscriptsize \xi}(t)}{dt}=\bm{\beta}(\bx_{\subscriptsize \xi}(t))
    \dfrac{d\bx_{\xi}(t)}{dt}=\bm{\beta}(\bx_{\xi}(t))
\end{equation}
with initial condition $\bx_{\scriptscriptstyle \xi}(0)=\bxi$, where $\bxi$ is a point on the inflow boundary $\Gamma_-$. %Each characteristic curve $\bx(t)$ depends on the starting point $\bxi\in \Gamma_-$ and may be conveniently denoted by $\bx_{\subscriptsize \xi} (t)$. 
Now, the PDE in \cref{pde1} becomes the following ordinary differential equations 
\begin{equation}\label{l-chara2}
    \dfrac{d\,u(\bx_{\scriptscriptstyle \xi}(t))}{dt} +\gamma u(\bx_{\scriptscriptstyle \xi}(t))=f(\bx_{\scriptscriptstyle \xi}(t))
\end{equation}
with initial condition $u(\bx_{\scriptscriptstyle \xi}(0))=g(\bxi)$ for all $\bxi\in \Gamma_-$. Denote by
\[
\Gamma_-^d=\left\{\bxi\in \Gamma_-: \, g \mbox{ is discontinuous at } \bxi\right\}
\]
the set of discontinuity of $g$. Then the discontinuity interface of the solution is given by
\[
\cI = \left\{\bx_{\scriptscriptstyle \xi}(t)\in \Omega: \, \bxi\in \Gamma_-^{d}, t\in [0,\infty), \mbox{ and } \bx_{\scriptscriptstyle \xi}(t) \mbox{ is the solution of \cref{l-chara}}\right\}.
\]
The interface $\cI$ partitions the domain $\Omega$ into submains $\{\Omega_i\}$, and one may compute the solution in each subdomain. This approach requires calculation of the discontinuity interface and is not applicable to nonlinear HCLs. 

To circumvent this difficulty, let us describe the approach developed in \cite{Cai2021linear}. To this end, introduce the directional differential operator 
%define the directional derivative of $u$ 
along the direction ${\bm\beta}$ by 
\begin{equation}\label{dd-d}
 D_{\scriptstyle{\bm\beta}} 
  \,v(\bx) = \lim_{\tau\to 0} \frac{u(\bx)-u\big(\bx - \tau {\bm{\beta}}(\bx)\big)}{\tau}
\end{equation}
for function $v$ whose directional derivative along $\bm\beta$ exists. %, where $\bar{\bm\beta}(\mathbf{x})=\bm\beta(\mathbf{x})/{|\bm\beta(\mathbf{x})|}$ is the unit vector along $\bm\beta(\mathbf{x})$ and $|{\bm\beta(\mathbf{x})}|$ is the magnitude of $\bm\beta(\mathbf{x})$. 
Then the linear advection-reaction equation with discontinuous solution may be written in the entire domain $\Omega$ as follows
\begin{equation}\label{pde1b}
    \left\{\begin{array}{rccl}
    %\bm{\beta}\cdot \!\nabla u 
    D_{\scriptstyle{\bm\beta}} \,u  + \gamma u &= & f &\text{ in }\,\, \Omega, \\[2mm]
    u&=&g &\text{ on }\,\, \Gamma_{-}. 
    \end{array}\right.
\end{equation}

Now, we are ready to introduce LS formulations of the linear advection-reaction problem based on \cref{pde1b}. To do so, denote the solution space by 
\begin{equation}\label{Vbeta}
   V_{\scriptstyle{\bm\beta}}= \left\{v\in L^2(\Omega): D_{\scriptstyle{\bm\beta}} \,v\in L^2(\Omega)\right\}, 
\end{equation}
%equipped with the norm \[ \|v\|_{\bm{\beta}}=\left(\|v\|_{0,\Omega}^2 + \|v_{\bm\beta} \|_{0,\Omega}^2\right)^{1/2}. \]
and define the following least-squares functional 
 \begin{equation}\label{ls1}
    \mathcal{L}(v;{\bf g}) = \|D_{\scriptstyle{\bm\beta}} \,v +\gamma\, v-f\|_{0,\Omega}^2 +  \|v-g\|_{-\bm\beta}^2, \quad \forall\,\, v\in V_{\bm\beta},
\end{equation}
where ${\bf g} = (f,g)$, and $\|\cdot\|_{0,\Omega}$ and $\|\cdot\|_{-\bm\beta}$ are the $L^2(\Omega)$ and the weighted $L^2(\Gamma_{-})$ norms on the domain and the inflow boundary given by
\[
\|v\|_{0,\Omega} = (v,v)^{1/2}=\left( \int_{\Omega}  v^2\,d\bx\right)^{1/2} \quad\mbox{and}\quad
\|v\|_{-\bm{\beta}} 
=\left<v,v\right>^{1/2}_{-\bm{\beta}} 
=\left( \int_{\Gamma_-} |\bm{\beta}\! \cdot \!\bm{n}|\, v^2\,ds\right)^{1/2},
\]
respectively.
%where $\hat{\bm{\beta}}$ is the unit vector in the direction $\bm{\beta}$.
Then the least-squares formulation of problem (\ref{pde1b}) studied in \cite{BoCh:2001, de2004least, BoGu16} is to seek $u\in V_{\bm\beta}$ such that
\begin{equation}\label{minimization1}
    \mathcal{L}(u;{\bf g}) = \min_{v\in V_{\bm\beta}} \mathcal{L}(v;{\bf g}).
\end{equation}
The well-posedness of the least-squares formulation in \cref{minimization1} was established in \cite{de2004least}.

%\textcolor{red}{
Problem \cref{minimization1} enforces the inflow boundary condition via a penalization in the weighted $L^2(\Gamma_-)$ norm. Alternatively, the condition can be imposed directly within the solution set, or both in the functional and the solution set. Specifically, define 
\begin{equation}\label{hat-ar}
   %\hat{\mathcal{L}}(v)= \|v_{\bm\beta} +\gamma\, v-f\|_{0,\Omega}^2 \quad\mbox{and}\quad
    {V}_{\scriptstyle{\bm\beta}}(g) = \left\{v\in {V}_{\scriptstyle{\bm\beta}} : v|_{\Gamma_-}=g\right\} = \left\{v\in L^2(\Omega): D_{\scriptstyle{\bm\beta}} \,v\in L^2(\Omega),\,v|_{\Gamma_-}=g\right\} .
\end{equation}
An equivalent least-square formulation is to find $u\in {V}_{\bm\beta}(g)$ such that
\begin{equation}\label{minimization11}
    {\mathcal{L}}(u;{\bf g}) = \min_{v\in {V}_{\bm\beta}(g)} {\mathcal{L}}(v;{\bf g}) \quad\mbox{or}\quad \hat{\mathcal{L}}(u) = \min_{ v\in {V}_{\bm\beta}(g)} \hat{\mathcal{L}}(v).
\end{equation}
where $\hat{\mathcal{L}}(v)= \|D_{\scriptstyle{\bm\beta}} \, v +\gamma\, v-f\|_{0,\Omega}^2$ is the first term of the LS functional ${\mathcal{L}}(v;{\bf g})$ in \cref{ls1}. Although the LS formulations in \cref{minimization1} and \cref{minimization11} are theoretically equivalent, practical considerations, 
especially the limited understanding of optimizers for non-convex problems, often motivate the use of first formulation in \cref{minimization11}, with an appropriately weighted second term in ${\mathcal{L}}(v;{\bf g})$.

\begin{remark}
The advection-reaction equation is often expressed in a conservative form as follows
\begin{equation}\label{pde1c}
    \left\{\begin{array}{rccl}
     \div ({\bm\beta} u) + \gamma u &= & f &\text{ in }\,\, \Omega, \\[2mm]
    u&=&g &\text{ on }\,\, \Gamma_{-}. 
    \end{array}\right.
\end{equation}
When the solution $u$ is discontinuous, the divergence operator $\div\!$ must be interpreted in a weak sense, as similarly defined in \cref{div}. We can then apply the LS principle to \cref{pde1c}
{\em (}see the subsequent section for details{\em )}.
\end{remark}

\subsection{Scalar Nonlinear Hyperbolic Conservation Laws}\label{s:hcls}

Let ${\bff}(u)=(f_1(u),...,f_d(u))$ be the spatial flux vector field, ${\Gamma}_-$ be the part of the boundary $\partial {\Omega} \times I$ where the characteristic curves enter the domain ${\Omega}  \times I \subset \R^{d+1}$, and the boundary data ${g}$ and the initial data $u_0$ be given scalar-valued functions defined on $\Gamma_-$ and $\Omega$, respectively. Consider the following scalar nonlinear hyperbolic conservation law 
\begin{equation} \label{pde2}
    u_t(\bx,t) + \sum\limits_{i=1}^d \dfrac{\partial f_i\left(u(\bx,t)\right)}{\partial x_i} = 0 \quad\text{in }\, {\Omega}  \times I
\end{equation}
with the inflow and initial conditions
\begin{equation} \label{BI2}
u=g \quad\mbox{on }\, {\Gamma}_{-} \quad\mbox{and}\quad u(\bx,0) = u_0(\bx) \quad\mbox{on }\, {\Omega},
\end{equation}
respectively, where $u_t$ is the partial derivative of $u$ with respect to the temporal variable $t$.
Without loss of generality, assume that $f_i(u)$ is twice differentiable for all $i\in \{1,\ldots, d\}$.

The solution of \cref{pde2} is often discontinuous due to a discontinuous initial or inflow boundary condition or shock formation. Hence, the strong form in \cref{pde2} is only valid where the solution is differentiable. Let $\cI$ denote the discontinuity interface of the solution. To fully characterize the behavior across $\cI$, an additional equation, the so-called Rankine-Hugoniot (RH) jump condition (see, e.g., \cite{leveque1992numerical, godlewski2013numerical}), is required:
\begin{equation}\label{RH}
  \left(\bff(u^+), u^+\right)\cdot \bn^+\big|_\cI + \left(\bff(u^-), u^-\right)\cdot \bn^-\big|_\cI=0,
\end{equation}
where $\bn^+=(\bn_x^+,n_t^+)$ and $\bn^-=(\bn_x^-,n_t^-)$ are the unit vector normal to the interface $\cI$ with opposite directions, and $u(\bx, t)$ have two different values on $\cI$, $u^+(\bx, t)$ and $u^-(\bx, t)$, defined by
\[
u^\pm(\bx, t)= \lim_{\epsilon \to 0} u(\bx-\epsilon\bn^\pm_x, t-\epsilon n^\pm_t) \quad\forall\,\, (\bx,t)\in \cI.
\]
Thus, a complete formulation of the scalar nonlinear hyperbolic conservation law must include both the PDE and the RH condition, yielding:
\begin{equation}\label{pde2a}
    \left\{\begin{array}{lcll}
    u_t(\bx,t) + \sum\limits_{i=1}^d \dfrac{\partial f_i\left(u(\bx,t)\right)}{\partial x_i} &= & 0, &\text{ in }\, \left({\Omega}  \times I\right)\setminus \cI, \\[4mm]
    \left(\bff(u^+), u^+\right)\cdot \bn^+ + \left(\bff(u^-), u^-\right)\cdot \bn^-&=& 0, &\text{ on }\,\cI ,
    \end{array}\right.
\end{equation}
supplemented with the inflow and initial conditions in \cref{BI2}. 

In practice, however, the interface $\cI$ is unknown {\it a priori}, making it extremely challenging to enforce the RH condition directly in numerical simulations. Traditional methods often rely on entropy conditions or artificial viscosity to approximate physically admissible solutions, but these approaches may introduce numerical artifacts such as smearing or oscillations near discontinuities. 

%To overcome the second difficulty, The concept of weak solutions of \cref{pde2} was introduced. $u\in L^\infty(\Omega\times I)$ is called a weak solution of (\ref{pde2}) if and only if\begin{equation}\label{weak}-(\bF (u), \nabla \varphi)_{0,\Omega\times I} + (\bn\cdot \bF (g),\varphi)_{0, \Gamma_-} + + (\bn\cdot \bF (u_0),\varphi)_{0, \Omega\times \{0\} } = 0,\quad\forall\,\, \varphi \in C^1_{\Gamma_+}(\Omega\times I).\end{equation}where $\bF(u)=(\bff(u),u)=(f_1(u),\ldots, f_d(u),u)$, $\Gamma_+=\partial \left(\Omega\times I\right) \setminus \left(\left(\Omega\times {0}\right) \cup \Gamma_-\right)$ is the outflow boundary and \[C^1_{\Gamma_+}(\Omega\times I)=\{\varphi\in C^1(\Omega\times I) :\, \varphi=0 \mbox{ on } {\Gamma_+}\}.\]
%One may deal with this difficulty through the weak solution defined in a variational form \cite{GoRa:96}. Since the collection of neural network functions is not a vector space, it is then difficult to discretize the weak formulation. Instead, 
To deal with this difficulty, let us introduce
the {\bf total flux}
\begin{equation}\label{F}
    \bF(u)=(\bff(u),u)=(f_1(u),\ldots, f_d(u),u).
\end{equation}
Denote by $\div\!$ the space-time divergence operator. When $u$ is differentiable at $(\bx,t)$, the classical definition of the divergence gives:
\begin{equation}\label{div-d}
\div\bF(u(\bx,t)) = u_t(\bx,t) + \sum\limits_{i=1}^d \dfrac{\partial f_i\left(u(\bx,t)\right)}{\partial x_i}.
\end{equation}
To interpret the \cref{RH} condition, we notice that it is indeed expresses a jump condition in the solution. From the perspective of of continuum mechanics, however, it represents the {\it continuity of the normal component of the total flux $\bF(u)$} across the interface $\cI$, i.e., 
\begin{equation}\label{RH2}
    \jump{\bF(u)\cdot\bn}_\cI \equiv \left(\bff(u^+), u^+\right)\cdot \bn^+\big|_\cI + \left(\bff(u^-), u^-\right)\cdot \bn^-\big|_\cI  = 0, 
\end{equation}
where $\jump{\cdot}_\cI$ denotes the jump over the interface $\cI$ and $\bn(\bx,t)$ is the unit normal vector to the interface. This observation motivates a generalized, weak definition of the divergence operator, applicable even at points where $u$ is discontinuous.
% This observation motivates us to define the divergence operator in a weak sense. 
Specifically, for any $(\bx,t)\in \cI$, 
% where $u(\bx,t)$ is discontinuous, 
we define the space-time divergence operator via the Gauss divergence theorem as:
\begin{equation}\label{div}
    \div\bF(u(\bx,t)) = \lim_{\epsilon \to 0} \dfrac{1}{|B_\epsilon(\bx,t)|} \int_{\partial B_\epsilon(\bx,t)} \bF(u)\cdot \bn \,dS,
\end{equation}
where $B_\epsilon(\bx,t) \in \R^{d+1}$ is a ball of radius $\epsilon$  centered at $(\bx,t)$, $\partial B_\epsilon(\bx,t)$ is the boundary of $B_\epsilon(\bx,t)$, and $\bn$ is the outward unit normal to $\partial B_\epsilon(\bx,t)$. 

\begin{theorem}\label{l:divS}
Let $u$ be a solution of \cref{pde2a} and \cref{BI2}, then the divergence of the total flux $\bF(u)$ vanishes on $\cI$, i.e.,
\begin{equation}\label{divS}
   \div {\bF} (u) =0\,\, \mbox{ in }\, \cI.
\end{equation}
\end{theorem}
\begin{proof}
For any $(\bx,t)\in \cI$, let $B_\epsilon(\bx,t)$ be a $\epsilon$-ball in $\Omega\times I$ centered at $(\bx,t)$. Then the interface $\cI$ partitions $B_\epsilon(\bx,t)$ into two subdomains denoted by $B^+_\epsilon(\bx,t)$ and $B^-_\epsilon(\bx,t)$ sharing part of $\cI$. It follows from \cref{RH2}, Gauss' divergence theorem, and the first equation of \cref{pde2a} that 
\begin{eqnarray*}
    \int_{\partial B_\epsilon(\bx,t)} \bF(u)\cdot \bn \,dS &=&  \int_{\partial B^+_\epsilon(\bx,t)} \bF(u)\cdot \bn \,dS + \int_{\partial B^-_\epsilon(\bx,t)} \bF(u)\cdot \bn \,dS \\[2mm]
    &=& \int_{B^+_\epsilon(\bx,t)} \div \bF(u) \,dS + \int_{B^-_\epsilon(\bx,t)} \div \bF(u) \,dS =0,
\end{eqnarray*}
which, together with the weak definition of the divergence operator in \cref{div}, implies \cref{divS}. This completes the proof of the lemma.
\end{proof}

%If $u$ is discontinuous at $(\bx,t)$, then $\div\bF(u(\bx,t))$ defined in \cref{div} leads to the continuity condition of the normal component of the space-time flux $\bF(u)$ that is identical to the RH jump condition. 
%\[\nabla\cdot {\bF} (u) = 0, \,\,\text{ in }\,\, ({\mathcal{R}^d}\times I) \setminus \textcolor{red}{\Gamma} \quad\mbox{and}\quad \jump{\bF(u)\cdot\bn}_{\textcolor{red}{\Gamma}} =0\]
%\[V_{{\bf f}}(u_0)=\{v\in L^2(\mathcal{R}^d\times I) :\, {\bf F}(v) \in H(\mbox{\textcolor{red}{div}};\mathcal{R}^d\times I), \, v(\bx,0)=u_0(\bx)\}\] 
%\[{\cal L}(v;\bff)=\|\mbox{\textcolor{red}{div}} \,{\bf F}(v)\|^2_{0, \mathcal{R}^d\times I}\]
%\[\Longrightarrow {\bf F}(v) \in H(\mbox{\textcolor{red}{div}};\mathcal{R}^d\times I)\]
%Find $u \in V_{{\bf f}}(u_0)$ such that \,\,\quad $u=\argmin\limits_{v\in V_{{\bf f}}(u_0)} {\cal L}(v;\bff)$

Using \cref{l:divS} and \cref{pde2a}, the nonlinear scalar hyperbolic conservation law has the following simplified form
\begin{equation} \label{pde2b}
   %\qquad\Longrightarrow \qquad 
   \div {\bF} (u) = 0 \,\,\text{ in }\, {\Omega}\times I\in \R^{d+1},
\end{equation}
supplemented by the inflow and initial conditions in \cref{BI2}.

As the linear case studied in \Cref{s:arp}, the inflow and initial conditions may be enforced either in the functional as penalization terms or in the solution set or in both the functional and the solution set. To this end,
denote the collection of square integrable vector fields whose divergence is also square integrable by
\[
H(\mbox{div};\Omega\times I)= \left\{\btau\in L^2(\Omega\times I)^{d+1} :\, \div\btau \in L^2(\Omega\times I)\right\}.
\]
Let us introduce the following two solution sets of (\ref{pde2b})
 \begin{eqnarray}\label{space}
 \left\{\begin{array}{ll}
     \cV_{\bff} = \left\{v\in L^2(\Omega\times I) :\, \bF(v)=(\bff(v),v)\in H(\mbox{div};\Omega\times I)\right\} \quad\mbox{ and }\, & \\[2mm]  
     {\cV}_{\bff}({\bf g})=\left\{v\in \cV_{\bff} :\, v|_{\Gamma_-}=g,\, v|_{\Omega\times \{0\}} =u_0\right\}, &
     \end{array}\right.
 \end{eqnarray}
 where ${\bf g}= (g, u_0)$.
Define two least-squares (LS) functionals as 
 \begin{equation}\label{ls2}
   \hat{\mathcal{L}} (v)=\| \div \bF(v)\|_{0,\Omega\times I}^2 \quad\mbox{and}\quad \mathcal{L}(v;{\bf g}) %\mathcal{L}(v;{g, u_0}) 
   =\hat{\mathcal{L}} (v) +  \|v-g\|_{0, \Gamma_-}^2 +  \|v-u_0\|_{0, \Omega \times \{0\}}^2,
\end{equation}
where $\|\cdot\|_{0,S}$ denotes the standard $L^2(S)$ norm for $S=\Omega\times I$, $\Gamma_-$, or $\Omega \times \{0\}$.
Now, the corresponding least-squares formulation is to either (1) seek $u\in \cV_{\bff}$ such that
\begin{equation}\label{minimization2a}
    \mathcal{L}(u;{\bf g}) = \min_{v\in \cV_{\bff}} \mathcal{L}(v;{\bf g})
\end{equation}
or (2) seek $u\in {\cV}_{\bff}(\bg)$ such that
\begin{equation}\label{minimization2b}
    {\mathcal{L}}(u;{\bf g}) = \min_{\small v\in {\cV}_{\bff}(\bf g)} {\mathcal{L}}(v;{\bf g}) \quad\mbox{or}\quad \hat{\mathcal{L}}(u) = \min_{\small v\in {\cV}_{\bff}(\bf g)} \hat{\mathcal{L}}(v).
\end{equation}
Comments on theory and practice of these three LS minimization problems in the previous section are valid here as well.

%\begin{proposition}Assume that $u\in L^\infty(\Omega\times I)$ is a piecewise $C^1$ function. Then $u$ is a weak solution of \cref{pde2b} if and only if $u$ is a solution of the minimization problem in either \cref{minimization2a} or \cref{minimization2b}.\end{proposition}\begin{proof}The proposition is a direct consequence of Theorem 2.5 in \cite{de2005numerical}.\end{proof}

\section{ReLU Neural Network and its Approximation to Discontinuous Functions}\label{s:NN}

%A neural network defines a new class of approximating functions which is suitable for approximating discontinuous functions with unknown interface locations such as solutions of nonlinear hyperbolic conservation laws. 

This section describes $l$-hidden-layer ReLU neural network as a set of continuous piecewise linear functions and illustrates its striking approximation power to discontinuous functions with {\it unknown} interface locations \cite{Cai2021linear, CCL2024}. 

%This chapter is restricted to one dimensional output $n_{l+1}=1$ for simplicity of presentation. Extension of materials covered by this chapter to multi-dimensional output $n_{l+1}>1$ is straightforward.

ReLU refers to the rectified linear activation function defined by
\begin{equation}\label{relu}
    \sigma(s) = \max\{0, s\} = \left\{\begin{array}{ll}
     s, & s>0,\\[2mm]
     0, & s\leq 0.
     \end{array}\right.
 \end{equation}
The $\sigma(s)$ is a continuous piecewise linear function with one {\it breaking} point $s=0$.
For $k=1,\ldots, l$, let $n_k$ denote the number of neurons at the $k^{th}$ hidden layer; denote by 
\[
\bb^{(k)}\in \R^{n_{k}}
\quad\mbox{and}\quad 
\bomega^{(k)} \in \R^{n_{k}\times n_{k-1}} 
\]
the biases and weights of neurons at the $k^{th}$ hidden layer, respectively. Their $i^{th}$ rows are denoted by $b_i^{(k)}\in \R$ and $\bomega^{(k)}_i\in \R^{n_{k-1}}$, that are the bias and weights of the $i^{th}$ neuron at the $k^{th}$ hidden layer, respectively, where $n_0=d$ or $d+1$ for applications in \Cref{s:problem}. Let 
\[
\bx^{(0)}=\bx \in \R^{d} 
\quad\mbox{or}\quad \bx^{(0)}=(\bx,t)\in \R^{d+1}.
\]
For $k=1,\ldots, l-1$, define vector-valued functions $\bx^{(k)}:\, \R^{n_{0}} \to \R^{n_k}$ by
\begin{equation}\label{layerdef}
  \bx^{(k)} = \bx^{(k)}\left(\bx^{(0)}\right) 
  = \sigma \left(\bomega^{(k)}\bx^{(k-1)}+\bb^{(k)}\right), % = \sigma \left(\bomega^{(k)}\bx^{(k-1)}\left(\bx^{(0)}\right)+\bb^{(k)}\right),
  %\quad\mbox{for } \bx^{(k-1)}\in \R^{n_{k-1}},
\end{equation}
where the application of the activation function $\sigma$ to a vector-valued function is defined component-wise.

A ReLU neural network with $l$ hidden layers and $n_k$ neurons at the $k^{th}$ hidden layer can be defined as the collection of continuous piecewise linear functions: 
\begin{equation}\label{network}  
\cM(l)=\left\{ c_0+\sum\limits_{i=1}^{n_l}c_i \sigma\left(\bomega^{(l)}_i\bx^{(l-1)}+b^{(l)}_i\right)  :  c_i,\, b^{(l)}_i\in \R, \,\, \bomega^{(l)}_i 
\in \cS^{n_{l-1}}\right\},
\end{equation}
where $\bx^{(l-1)}\left(\bx^{(0)}\right)$ is recursively defined as in \cref{layerdef}, and $\cS^{n_{l-1}}$ denotes the unit sphere in $\R^{n_{l-1}}$. The constraint that the weight vector of each neuron lies on the unit sphere arises from normalization for the ReLU activation function (see \cite{LiuCai1}). This restriction can narrow the set of solutions for a given approximating problem.
The total number of parameters of $\cM(l)$ is given by
\begin{equation}\label{DoF}
   M(l) =(n_l+1) + \sum^{l}_{k=1} n_{k}\times (n_{k-1}+1) .
\end{equation}

We now explain why any function $v\left(\bx^{(0)}\right)$ in $\cM(l)$ is always a continuous piecewise linear function, regardless of the dimension, the number of layers, or the number of neurons per layer. %Without loss of generality, let $\bx^{(0)}=\bx \in \R^{d}$. 
To illustrate this, we first consider the shallow NN case, i.e., $l=1$. For simplicity, let $\bx^{(0)}=\bx \in \R^{d}$, then any NN function $v^{(1)}\left(\bx\right)\in \cM(1)$ has the form of
\begin{equation}\label{v1}
v^{(1)}\left(\bx\right)= c_0+\sum\limits_{i=1}^{n_1}c_i \sigma\left(\bomega^{(1)}_i\bx+b^{(1)}_i\right).
\end{equation}
Since the ReLU activation function $\sigma(s)$ is itself a continuous piecewise linear function with a single {\it breaking point} at $s=0$, each neuron $\sigma\left(\bomega^{(1)}_i\bx+b^{(1)}_i\right)$ is a continuous piecewise linear function with a {\it breaking hyperplane} (see \cite{Cai2021linear,LiuCai1}):
\begin{equation}\label{P-i}
\mathcal{P}_i^{(1)}=\left\{\bx\in \Omega\subset \R^d:\, {\bomega}^{(1)}_{i}\cdot\bx + b^{(1)}_i=0 \right\}.
\end{equation}
These hyperplanes, together with the boundary of the domain $\Omega$, generate a {\it physical partition} ${\cal K}^{(1)}$ of $\Omega$ as described in \cite{LiuCai1, Cai2023AI}. This partition ${\cal K}^{(1)}$ is uniquely determined by the parameters of the first layer $\left\{ \left( \bomega_i^{(1)}, b_i^{(1)} \right) \right\}_{i=1}^{n_1}$ and consists of irregular, polygonal subdomains of $\Omega$ (see \cref{breaking5} below for a visual example for such partition). The NN function $v^{(1)}(\bx)$ in \cref{v1} is then a continuous piecewise linear function with respect to ${\cal K}^{(1)}$.  

Any NN function $v^{(2)}\left(\bx\right)$ in $\cM(2)$ has the form
\begin{equation}\label{v2}
v^{(2)}\left(\bx\right)= c_0+\sum\limits_{i=1}^{n_2}c_i \sigma\left(\bomega^{(2)}_i\bx^{(1)}+b^{(2)}_i\right),
\end{equation}
where $\bx^{(1)} = \sigma \left(\bomega^{(1)}\bx +\bb^{(0)}\right)\in \cM(1)^{n_1}$ is a vector-valued function whose components are continuous piecewise linear functions. Hence, each pre-activated NN function 
\[
g_i^{(2)}(\bx) = \bomega^{(2)}_i\bx^{(1)}+b^{(2)}_i, 
\]
being an affine function of $\bx^{(1)}$, is also a continuous piecewise linear function of $\bx$, defined with respect to the physical partition ${\cal K}^{(1)}$. The zero level set of $g_i^{(2)}(\bx)$, denoted by
\begin{equation}\label{P-i2}
\mathcal{P}_i^{(2)}=\left\{\bx\in \Omega\subset \R^d:\, \bomega^{(2)}_i\bx^{(1)}+b^{(2)}_i=0 \right\}
\end{equation}
defines a breaking poly-hyperplane that refines the partition ${\cal K}^{(1)}$ (see \cite{Cai2023AI} for details). It follows that $v^{(2)}\left(\bx\right)$ is a continuous piecewise linear function with respect to the refined physical partition ${\cal K}^{(2)}$, generated by incorporating the breaking poly-hyperplanes $\left\{ \cP^{(2)}_i \right\}_{i=1}^{n_2}$ into $\mathcal{K}^{(1)}$. By induction, any ReLU NN function $v^{(l)}\left(\bx\right)\in \cM(l)$ is a a continuous piecewise linear function with respect to the physical partition ${\cal K}^{(l)}$, which is determined by the network parameters across all hidden layers.

One of the key challenges in numerically solving HCLs lies in the fact that the location of the discontinuities in the solution is typically unknown in advance. To capture these unknown interfaces, traditional mesh-based numerical methods require fine meshes, whick leads to high computational cost and often produces spurious oscillations near the interface, an artifact commonly referred to as the Gibbs phenomenon. In the remainder of this section, we highlight a remarkable approximation property of ReLU NN functions for representing piecewise constant functions with unknown discontinuity interfaces.

To this end, let $\chi(\bx)$ be a piecewise constant function defined on $\Omega \in \R^{d}$, given by
\begin{equation}\label{pc}
   \chi(\bx)= \sum_{i=1}^k \alpha_i {\bm 1}_{\Omega_i}(\bx), 
\end{equation}
where $\left\{\Omega_i\right\}_{i=1}^k$ forms a partition of the domain $\Omega$ and ${\bm 1}_{\Omega_i}(\bx)$ denotes the indicator function of the subdomain $\Omega_i$. For simplicity of presentation, consider a special case that $k=2$, $\alpha_1=0$, and $\alpha_2=1$, i.e.,
\begin{equation}\label{step function}
 \chi(\bx)={\bm 1}_{\Omega_2}(\bx)=\left\{\begin{array}{rl}
 0, & \bx \in \Omega_1,\\[2mm]
 1, & \bx \in \Omega_2, 
 \end{array}
 \right.   
\end{equation}
where $\Omega_1$ and $\Omega_2$ are open, connected subdomains of $\Omega$ satisfying 
\[
\Omega_1\cap  \Omega_2=\emptyset 
\quad\mbox{and}\quad
\bar{\Omega}=\bar{\Omega}_1\cup  \bar{\Omega}_2.
\]
Let $\partial \Omega_i$ be the boundary of the subdomain $\Omega_i$, and assume the interface $\cI=\partial \Omega_1\cap \partial\Omega_2$ is a $C^0$ surface with finite $(d-1)$-dimensional measure, i.e., $|\cI|<\infty$. Although this case only considers a single interface, the extension to general piecewise constant functions as in \cref{pc} with multiple discontinuities is straightforward.
% (Extension to general piecewise constant function $\chi(\bx)$ defined in \cref{pc} is straightforward.)}
%Denote by $\Gamma=\partial \Omega_1\cap \partial\Omega_2$ the discontinuity interface of $\chi(\bx)$, where $\partial \Omega_i$ is the boundary of the subdomain $\Omega_i$. Assume that the interface $\Gamma$ is $C^0$ and that its $(d-1)$-dimensional measure $|\Gamma|$ is finite.

\begin{figure}[htbp]\label{1dstep}
\centering
\subfigure[$p_1(\bx)$\label{1dstep1}]{
\begin{minipage}[t]{0.4\linewidth}
\centering
\includegraphics[width=2.3in]{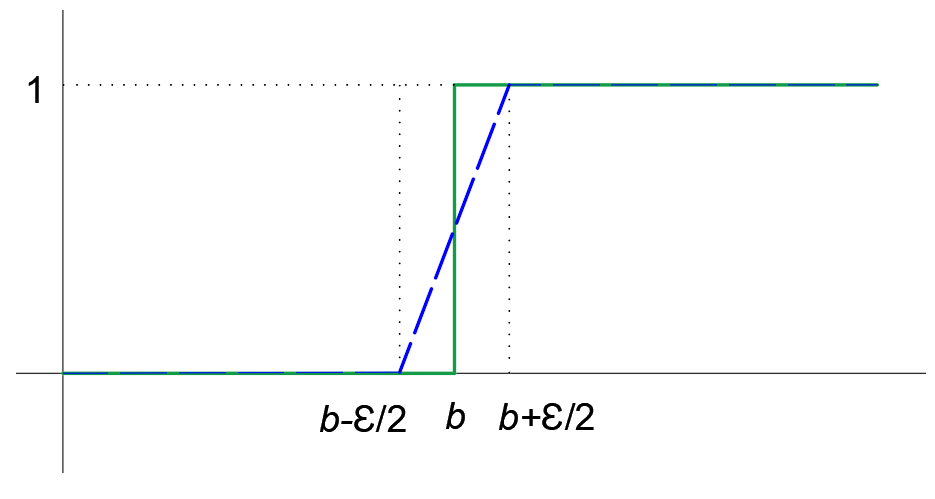}
\end{minipage}%
}%
\hspace{0.25in}
\subfigure[$p_2(\bx)$\label{1dstep2}]{
\begin{minipage}[t]{0.4\linewidth}
\centering
\includegraphics[width=2.3in]{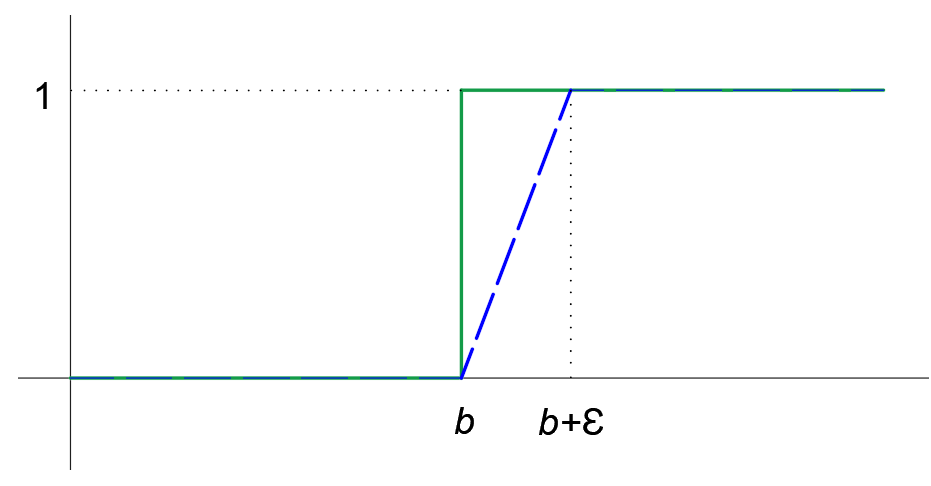}
\end{minipage}%
}%
\caption{NN Approximation of a 1D unit step function.}
\end{figure}

The complexity of the discontinuity interface $\cI$ plays a crucial role in determining both the accuracy of the approximation and the required architecture of NN. To gain a clear understanding, we begin by examining the simplest non-trivial case, when the interface is a hyperplane, and defer a discussion of more general interfaces to Remark~\ref{r:1}.
% we restrict our discussion to the simplest case: the hyperplane interface, and we then comment on the general case in \cref{r:1}. 
% To this end, 
%In this chapter, we assume that $\Gamma$ can be approximated by a connected series of hyperplanes with a prescribed accuracy $\varepsilon>0$ (see \cref{general} and \cite{CCL2024}).
%we first show by construction that ReLU NN functions can approximate step functions, whose interfaces are hyperplanes, within any prescribed accuracy $\varepsilon > 0$. One construction uses a two-layer NN with only two neurons and the other employs a three-layer NN with one neuron per hidden layer. We then comment on approximation results of ReLU NN functions to piecewise constant functions with slightly general interface. 
Specifically, consider the case where the interface $\cI$ is a portion of a hyperplane defined by
\[
\cI = \left\{\bx\in \Omega\subset\R^d :\, \mathbf{a}\cdot\mathbf{x}=b \right\},
\]
where $\mathbf{a}$ is the unit vector normal to $\cI$ pointing to the subdomain $\Omega_2$. To approximate the step function $\chi(\bx)$ associated with this hyperplane interface, we introduce two representative ReLU NN constructions, both capable of capturing the discontinuity sharply and efficiently.
For any given $\varepsilon >0$, define the following NN functions as in \cite{Cai2021linear} and \cite{CCL2024}:
\begin{equation}\label{Gamma HP}
   p_1(\bx)= \dfrac{1}{2\varepsilon} \big(\sigma(\mathbf{a}\!\cdot\!\mathbf{x}-b+\varepsilon) - \sigma(\mathbf{a}\!\cdot\!\mathbf{x}-b-\varepsilon) \big)  \quad\mbox{or} \quad p_2(\mathbf{x})=1-\sigma\left(-\frac{1}{\varepsilon}\sigma(\mathbf{a}\!\cdot\!\mathbf{x}-b)+1\right).
\end{equation}
Here, $p_1(\bx)$ is achieved by a two-layer ReLU NN with just two neurons, while $p_2(\bx)$ corresponds to a three-layer NN with one neuron per hidden layer. In one dimension, they are illustrated in \cref{1dstep}. Define the narrow transition regions as
\[
\Omega_{p_1}=\left\{ \bx\in\Omega: \,
    -\varepsilon<\mathbf{a}\cdot \mathbf{x}-b<\varepsilon \right\} 
    \quad\mbox{and}\quad
 \Omega_{p_2}=\left\{ \bx\in\Omega: \,
    0<\mathbf{a}\cdot \mathbf{x}-b<\varepsilon \right\}.    
\]
Then, both $p_1(\bx)$ and $p_2(\bx)$ serve as smoothed approximations to the ideal step function $\chi(\bx)$, with the following piecewise-defined formula:
\begin{equation}\label{Gamma HP2}
 p_1(\bx)=\left\{\begin{array}{cc}
    0, & \mbox{in }\, \Omega_1\setminus \Omega_{p_1}, \\[3mm]
    \dfrac{\mathbf{a}\cdot \mathbf{x}-b+\varepsilon}{2\varepsilon}, & \mbox{in }\, \Omega_{p_1}, \\[3mm]
    1, & \mbox{in }\, \Omega_2\setminus \Omega_{p_1}
\end{array}\right. \,\mbox{ and }\, p_2(\bx)=\left\{\begin{array}{cc}
    0, & \mbox{in }\, \Omega_1\setminus \Omega_{p_2}, \\[3mm]
    (\mathbf{a}\cdot \mathbf{x}-b)/\varepsilon, & \mbox{in }\, \Omega_{p_2}, \\[3mm]
    1, & \mbox{in }\, \Omega_2\setminus \Omega_{p_2}.
\end{array}\right.   
\end{equation}

\begin{lemma}\label{l:I1}
There exists a positive constant $C$ such that for all $r\in (0,\infty)$, we have
\begin{equation}\label{error-bound}
     \|\chi-p_1\|_{L^r(\Omega)}\le C\,|\cI|^{1/r}\varepsilon^{1/r}
    \quad\mbox{and}\quad 
    \|\chi-p_2\|_{L^r(\Omega)}\le C\,|\cI|^{1/r}\varepsilon^{1/r},
\end{equation}
where $|\cI|$ is the $(d\!-\!1)$-dimensional measure of the interface $\cI$.
\end{lemma}
\begin{proof}
By \cref{step function} and \cref{Gamma HP2}, for $i=1,2$, we have
\[
\chi(\bx)-p_i(\mathbf{x})=0, \mbox{ in }\, \Omega\setminus\Omega_{p_i} \quad\mbox{and}\quad \big|\chi(\bx)-p_i(\mathbf{x})\big|^r \leq 1, \mbox{ in }\, \Omega_{p_i},
\]
%\begin{equation}\label{error1}    \big|\chi(\bx)-p_1(\mathbf{x})\big| \left\{\begin{array}{cc}    =0, & \mbox{in }\, \Omega\setminus \Omega_{p_1}, \\[2mm]    \leq 1, & \mbox{in }\, \Omega_{p_1}\end{array}\right. \,\mbox{ and }\,\big|\chi(\bx)-p_2(\mathbf{x})\big|\left\{\begin{array}{cc}    =0, & \mbox{in }\, \Omega\setminus \Omega_{p_2}, \\[2mm]    \leq 1, & \mbox{in }\, \Omega_{p_2},\end{array}\right. \end{equation}
%It is then easy to see that \begin{equation}\label{error2} \left|\chi(\mathbf{x})-p_1({\mathbf{x}})\right|^r \le 1, \,\,\forall \,\bx\in \Omega \setminus \Omega_{p_1} \quad\mbox{and}\quad \left|\chi(\mathbf{x})-p_2({\mathbf{x}})\right|^r \le 1, \,\,\forall \,\bx\in \Omega \setminus \Omega_{p_2}, \end{equation}
which, together with the facts that $\big|\Omega_{p_i}\big|\leq C\,\big|\cI\big|\varepsilon$, 
%\[\big|\Omega_{p_1}\big|\leq C\,\big|\Gamma\big|\varepsilon\quad\mbox{and}\quad\big|\Omega_{p_2}\big|\leq C\,\big|\Gamma\big|\varepsilon,\]
implies the validity of \cref{error-bound}. This completes the proof of the lemma.
\end{proof}

\begin{figure}[htbp]\label{general}
\centering
\subfigure[The interface $\Gamma$\label{general1}]{
\begin{minipage}[t]{0.4\linewidth}
\centering
\includegraphics[width=1.8in]{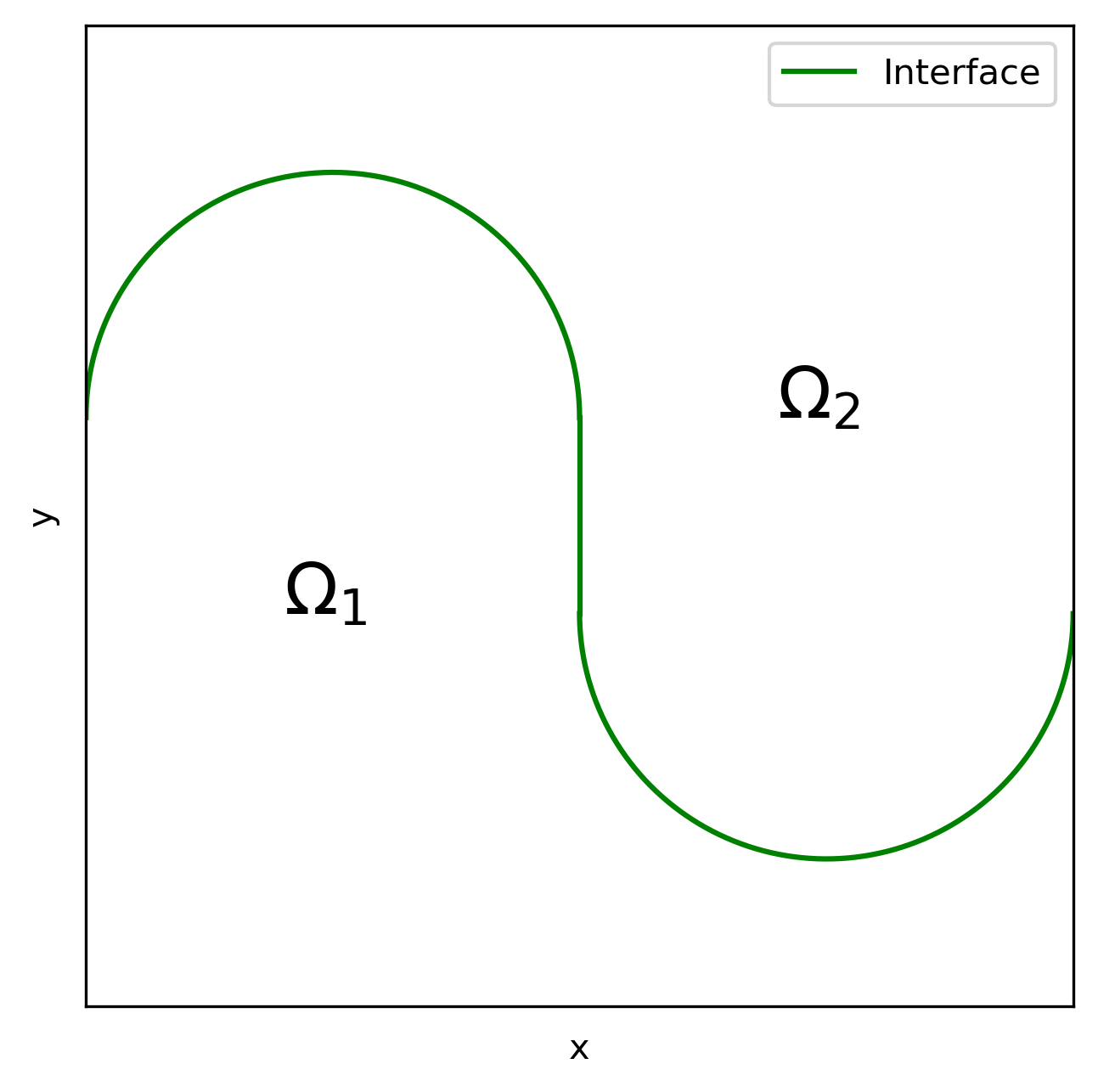}
\end{minipage}%
}%
\hspace{0.2in}
\subfigure[An approximation of the interface by connected series of hyperplanes\label{general22}]{
\begin{minipage}[t]{0.4\linewidth}
\centering
\includegraphics[width=1.8in]{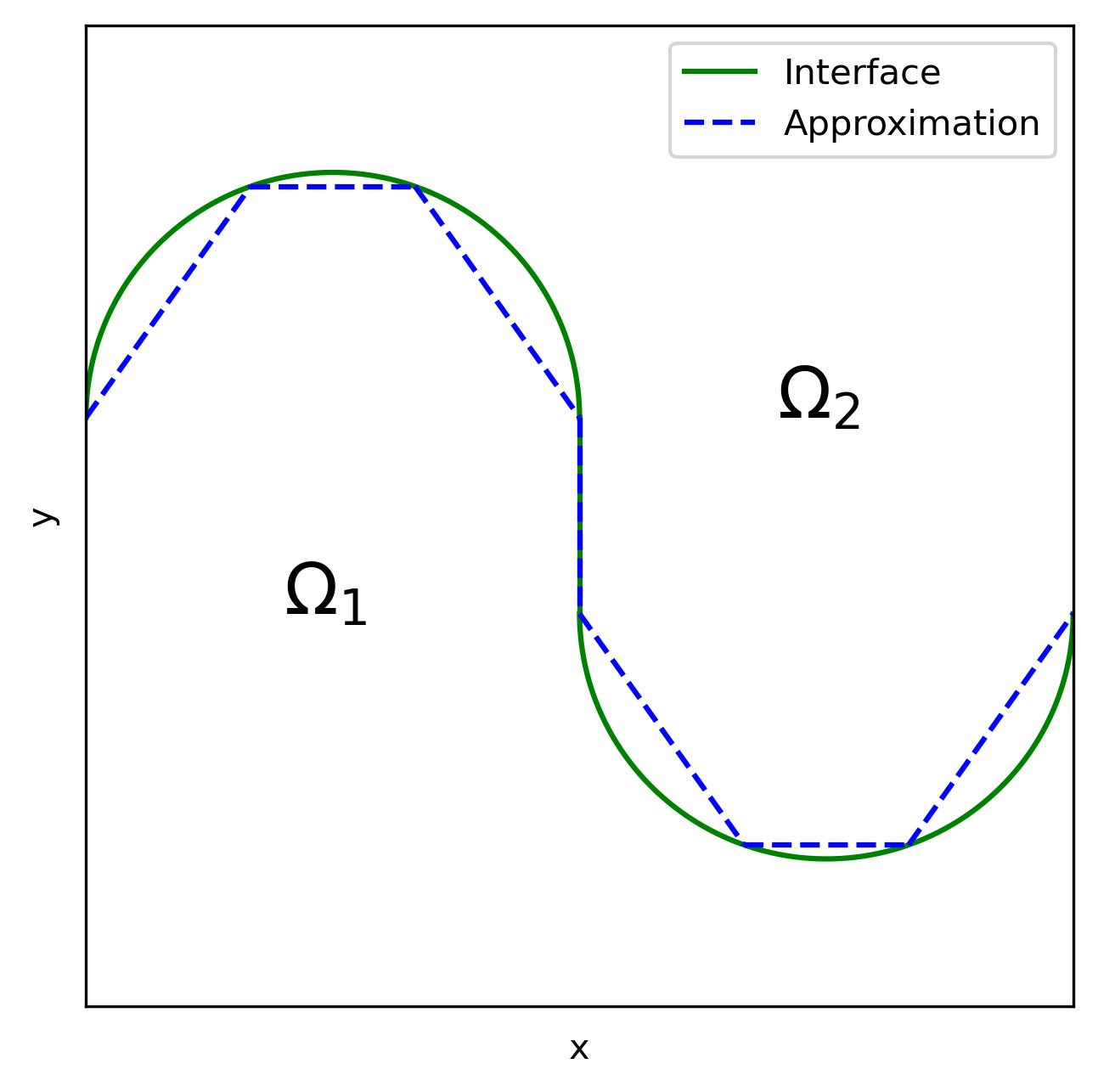}
\end{minipage}%
}%
\caption{Approximation of the interface $\Gamma$}
\end{figure}

These explicit constructions demonstrate the expressive power of ReLU NNs in approximating discontinuous functions, with simple architectures sufficing to localize and resolve interfaces aligned with hyperplanes. In the next remark, we discuss the extension of this approach to more general, non-planar interfaces.

\begin{remark}\label{r:1}
When the interface $\cI$ is not a hyperplane but can be approximated, to a prescribed accuracy $\varepsilon >0$, by a connected sequence of hyperplanes {\em(}see {\em \cref{general} and \cite{CCL2024})} in the maximum norm, %\textcolor{red}{in the maximum norm}, 
the indicator function $\chi(\bx) = {\bm 1}_{\Omega_2}(\bx)$ defined in \cref{step function} can still be effectively approximated by a ReLU NN of a given architecture, satisfying the error bound in \cref{error-bound}. 

More precisely, based on the two-layer ReLU NN approximation $p_1(\bx)$ in \cref{Gamma HP}, we showed in {\em \cite{Cai2023linear}} that a ReLU NN with at most $\lceil \log_2(d+1)\rceil+1$ layers suffices to achieve the desired approximation accuracy $\varepsilon$. However, {\em \cite{Cai2023linear}} does not provide an estimate on the minimum number of neurons at each layer. This result is obtained by constructing a continuous piecewise linear (CPWL) function with a sharp transition layer of $\varepsilon$ width, and combining the main findings of \cite{arora2016understanding} regarding the ReLU NN representation of CPWL functions.

Leveraging  the three-layer ReLU NN approximation $p_2(\bx)$ in \cref{Gamma HP}, it was shown in {\em \cite{CCL2024}} that $\chi(\bx)$ can be approximated with the comparable  accuracy using a three-layer ReLU NN. In this case, the number and locations of neurons at the first layer are determined by the hyperplanes used to approximate the interface, while the number of neurons in the second layer is dependent on the convexity of the interface {\em (}see {\em Theorem~3.2} in {\em \cite{CCL2024})}.
\end{remark}

%\begin{remark}\label{r:2} Let $\left\{\Omega_i\right\}_{i=1}^k$ be a partition of the domain $\Omega$. Let $\chi(\bx)$ be a piecewise constant function with respect to the partition with $\chi(\bx)=\alpha_i$ in $\Omega_i$ for $i=1,\ldots,k$. Then we have \[\chi(\bx)= \sum_{i=1}^k \alpha_i {\bm 1}_{\Omega_i}(\bx),\] where ${\bm 1}_{\Omega_i}(\bx)$ is the indicator function of the subdomain $\Omega_i$. As indicated in \cref{r:1}, each indicator function may be approximated by a ReLU function with a prescribed accuracy, and so is $\chi(\bx)$. \end{remark}

\section{Least-Squares Neural Network (LSNN) Method}\label{s:LSNN}

This section introduces the least-squares neural network (LSNN) method for solving advection-reaction equations in \cref{pde1} and scalar nonlinear hyperbolic conservation laws in \cref{pde2b} based on the equivalent least-squares formulations in \cref{minimization1}/\cref{minimization11} and \cref{minimization2a}/\cref{minimization2b}, respectively. To evaluate the least-squares functionals, we discuss efficient numerical integration in \Cref{s:NI} and physics-preserved numerical differentiation in \Cref{s:DDO}. Finally, the LSNN method is defined in \Cref{s:LSNNs}.

\subsection{Numerical Integration}\label{s:NI}
The evaluation of the least-squares functionals $\mathcal{L}(v;{\bf f})$ and $\hat{\mathcal{L}}(v;{\bf f})$, as defined in \cref{ls1} or \cref{ls2}, involves integrations over the computational domain $\Omega\subset \R^d$ or $\Omega\times I\subset \R^{d+1}$ ($d=1,2$, or $3$), as well as portions of their boundaries. In practice, these integrals must be approximated using numerical quadrature. This section outlines the basic setup for numerical integration and 
discusses key considerations for its implementation within the LSNN framework.

Let 
\[
{\cal T}=\{K\, :\, K\mbox{ is an open subdomain of } \Omega\}
\]
denote a partition of the domain $\Omega$, referred to as the {\em integration mesh}. %The partition satisfies
Here, the partition means that the union of all subdomains of ${\cal T}$ equals the whole domain $\Omega$ and that any two distinct subdomains of ${\cal T}$ have no intersection; more precisely,
\[ 
\bar{\Omega} = \cup_{K\in {\cal T}} \bar{K} \quad\mbox{and}\quad
 K\cap T = \emptyset, \quad \forall\,\, K,\, T \in {\cal T}.
\]
A composite quadrature over ${\cal T}$ is then written as
\[
\sum_{K\in\mathcal{T}} \mathcal{Q}_K(w)\approx \sum_{K\in\mathcal{T}} \int_K w(\bx)\,d\bx =\int_\Omega w(\bx)\,d\bx,
\]
where $\mathcal{Q}_K(w)\approx \int_K w(\bx)\,d\bx$ represents a quadrature rule over subdomain $K$. The specific quadrature rule $\mathcal{Q}_K$ may vary across different elements $K\in \mathcal{T}$, and can be chosen from standard formulas such as Gaussian quadrature or Newton--Cotes formulas, including the midpoint, trapezoidal, or Simpson’s rule (see \cite{Stroud71}). For instance, using the midpoint rule for all $K\in \mathcal{T}$, $\mathcal{Q}_K(w)=w(\mathbf{x}_K)\lvert K\rvert$, where $\mathbf{x}_K$ is the centroid of element $K$ and $\lvert K\rvert$ denotes its $d$- or $(d\!+\!1)$-dimensional measure. 
% First, $\mathcal{Q}_K$ may vary on $K\in \mathcal{T}$. Second, its choice is one of the standard quadrature rules like the Gaussian quadrature or Newton--Cotes formulas such as the midpoint, trapezoidal, or Simpson rule (see \cite{Stroud71}). In the case of the midpoint rule for all $K\in \mathcal{T}$, $\mathcal{Q}_K(w)=w(\mathbf{x}_K)\lvert K\rvert$, where $\mathbf{x}_K$ is the centroid of $K$ and $\lvert K\rvert$ is the $d$- or $(d\!+\!1)$-dimensional measure of $K$. 

In the context of LSNN, the integrands typically depend on NN approximations of the solution $u$ to the underlying PDE. These NN approximations are continuous piecewise linear functions defined over physical partitions (see \Cref{s:NN}), which are generally unknown in advance and evolve dynamically during training. Moreover, the true solution
$u$ itself is unknown and may exhibit localized features such as steep gradients, discontinuities, or singularities.

%\textcolor{red}{In the application of NN-based discretization, integrands depend on NN approximations to the solution $u$ of the underlying PDE. Each ReLU NN approximation is a continuous piecewise linear function with respect to a partition of the domain. Such a partition is referred to as the physical partition of the current NN approximation in \cite{LiuCai1} and is determined by the weights and biases of the neurons. Moreover, the physical partition is highly irregular (see \cref{breaking5} on page 20 and, e.g., Figures 7.4(b,d,f,g,h) and 7.5(d) in \cite{LiuCai1}) and moves during trainings. Additionally, the solution $u$ is unknown and has some local features. }

Because of these considerations, adaptive numerical integration was introduced in \cite{LiuCai1} (see Algorithm~5.2) and in \cite{LiCaRa23} (see Algorithm~3.1). Below, we briefly outline the adaptive mesh refinement strategy for numerical integration with a fixed NN in Algorithm~4.1 for problem \cref{pde1}. 

As usual, we begin with a uniform and coarse partition $\mathcal{T}^\prime$ of the domain $\Omega$. Based on this initial integration mesh, we apply the standard adaptive numerical quadrature procedure to refine it according to the given integrands $f$ and $g$. The goal is to construct an initial integration partition $\cT$ such that %$\int_\Omega f^2\, d\bx$ can be approximated accurately by the corresponding composite quadrature rule; specifically, for a prescribed accuracy $\varepsilon> 0$, we have
\[
\left| \int_\Omega f^2\, d\bx - \sum_{K\in \cT} \mathcal{Q}_K\left(f^2\right)\right| \leq \varepsilon \quad\mbox{and}\quad \left|\int_{\Gamma_-} g^2 dS  - \sum_{E\in {\cal E}_-} \mathcal{Q}_E\left(g^2\right)\right| \leq \varepsilon,
\]
%In a similar fashion, starting with $\cT^{\prime\prime}$, we create an initial  integration partition $\cT^{\prime\prime}\subset \cT$ such that $\int_{\Gamma_-} g^2 dS$ can be approximated accurately . 
where ${\cal E}_-$ denotes the collection of inflow boundary faces of $\cT$ (see \cref{II}).

%Assume that the inflow boundary data $g$ can be approximated with a prescribed accuracy by a continuous piecewise linear function with respect to the partition $\mathcal{T}$. 
Let $u_{_\mathcal{T}}$ be the NN approximation obtained using this initial integration mesh $\cT$. For each subdomain $K\in \mathcal{T}$, we define a local error indicator
\[
\eta_{_K}=\|D_{\scriptstyle{\bm\beta}, \tau} u_{_{\mathcal{T}}} +\gamma\, u_{_{\mathcal{T}}}-f\|_{0,K},
\]
where $D_{\scriptstyle{\bm\beta}, \tau}$ is a discrete directional derivative operator as defined in \cref{finite_diff}. The global error estimator is then given by $\eta=\left(\sum\limits_{K\in\cT} \eta_{_K}^2\right)^{1/2}$. These two error indicators guide the adaptive refinement process, which iteratively updates the integration mesh to better resolve regions with large local error contributions. The refinement procedure is summarized in Algorithm~4.1.
% The adaptive mesh refinement is summarized in Algorithm~4.1. %\cref{a:4.1}.

\begin{algorithm}\label{a:4.1}
{\bf {\sc \bf Algorithm 4.1}} Adaptive Quadrature Refinement (AQR) with a fixed NN.
\begin{itemize}
     \item[(1)] for each $K\in\cT$, compute the local error indicator $\eta_{K}$;
    \item[(2)] mark $\cT$ by the either bulk or average marking strategy (see, e.g., \cite{LiuCai1}) and refine marked subdomain to obtain a new partition $\cT^\prime$;
    \item[(3)] compute new NN approximation $u_{\mathcal{T}^\prime}$ on the refined integration mesh $\mathcal{T}^\prime$;
    \item[(4)] for a given parameter $\gamma \in (0, 1)$, if $\eta (u_{_{\cT^\prime}}) \leq \gamma \eta (u_{_{\cT}})$,
    then go to Step (1) with $\cT = \cT^\prime$; otherwise, output $\cT$.
\end{itemize}
\end{algorithm}

\smallskip

%As indicated in \cite{LiuCai1, LiCaRa23}, 
For a given parameter $\gamma \in (0, 1)$, the stopping criterion used in Step (4) of Algorithm~4.1 indicates that the adaptive procedure terminates when further refinement of the integration mesh produces a negligible reduction in residual. In other words, if the additional refinement does not significantly improve the accuracy of the global error evaluation for the current NN approximation, the algorithm halts and outputs the current integration partition.
% means that the adaptive procedure stops if further refinement of integration partition does not reduce much the size of the residual. % the current NN. %, the adaptive quadrature stops and outputs the current integration mesh.
    
\begin{remark}
In the case where computational cost is not an issue, one may use a uniform partition $\cT$ that is fine enough to approximate the unknown solution well using a piecewise polynomial.
\end{remark}

%\begin{remark}\label{int3} The LSNN approximation $u_N\in \cM(L,n)$ defined in \cref{L-NN} is CPWL with respect to a partition of the domain $\Omega$, referred to as the physical partition in \cite{LiuCai1, LiuCai2, Cai2021DeepAdaptive}. The partition $\mathcal{T}$ is a ``mesh'' for numerical integration and is completely different from the physical partition of $u_N$. Therefore, the partition $\mathcal{T}$ differs from meshes of traditional numerical methods. Nevertheless, the partition $\mathcal{T}$ and the corresponding quadrature $\mathcal{Q}_K$ are important for accuracy of the approximation $u_N$ by providing accurate information of the exact solution.  \end{remark}

\subsection{Physics-Preserved Numerical Differentiation}\label{s:DDO}

Solutions to the hyperbolic PDEs in \cref{pde1} and \cref{pde2} can exhibit discontinuities, making conventional numerical or auto-differentiations along coordinate directions, based directly on \cref{pde1} and \cref{pde2}, inadequate. This section introduces physics-preserved numerical differentiation techniques derived from the equivalent formulations in \cref{pde1b} and \cref{pde2b}, as proposed in \cite{Cai2021linear, Cai2023nonlinear}. 

When the solution $u$ is discontinuous, as discussed in \Cref{s:arp} and \Cref{s:hcls}, both the directional derivative $u_{\bm\beta}(\bx)$ and the divergence of the total flux $\div {\bF} (u)$ can still be defined via limit processes, as given in \cref{dd-d} and \cref{div}, respectively. Approximating these quantities through any reasonable limiting process leads to what we call {\em physics-preserved numerical differentiation}. 

Based on \cref{dd-d}, for any function $v(\mathbf{x})$ that admits a directional derivative in the $\bm\beta$ direction at $\bx\in \Omega$, we define the discrete directional differential operator $D_{\scriptstyle{\bm\beta},\tau}$ as
%Based on \cref{dd-d}, for any $\bx\in \Omega$, define the discrete differential operator $D_{\bm\beta}$ by
 \begin{equation}\label{finite_diff}
     D_{\scriptstyle{\bm\beta},\tau}\,v(\mathbf{x}) %\coloneqq 
     = \frac{v(\bx)-v\big(\bx - \tau {\bm{\beta}}(\bx)\big)}{\tau} %\approx v_{\bm\beta}(\bx),
\end{equation}
for $0<\tau \ll 1$. This operator provides an upwind finite difference approximation to the directional derivative in the direction of $\bm\beta$ with step size $\tau$. 
Intuitively, for sufficiently small $\tau$, this scheme avoids crossing the discontinuity interface, thereby preserving the physical meaning of the directional derivative. In contrast, standard numerical or automatic differentiation based on $\sum\limits_{i=1}^d\beta_i u_{x_i}$, where $u_{x_i}$ denotes the partial derivative of $u$ with respect to $x_i$, involves values from both sides of the interface and thus fails to capture the correct behavior at discontinuities.
This can be illustrated by a simple model problem in $\Omega=(0,1)\times (0,1)\subset\R^2$ with $\bm\beta=(1,1)^T$. If the inflow boundary condition $g$ is discontinuous at the corner point $(0,0)$, the resulting discontinuity interface is $\cI=\{(x,t)\in\Omega : x+t=0\}$. In this scenario, the standard differentiation $u_x+u_t$ necessarily evaluates across the interface $\cI$, leading to non-physical or spurious results.

% Intuitively, for sufficiently small $\tau$, it computes the derivative without crossing the discontinuity interface, while the standard numerical/auto differentiation based on $\sum\limits_{i=1}^d\beta_i u_{x_i}$ are computed using values on both sides of the interface, where $u_{x_i}$ denotes the partial derivative of $u$ with respect to $x_i$. This can be seen in the simplest model problem defined in $\Omega=(0,1)\times (0,1)\subset\R^2$ with $\bm\beta=(1,1)^T$. If the inflow boundary condition $g$ is discontinuous at $(0,0)$, then the discontinuity interface is $\Sigma=\{(x,t)\in\Omega : x+t=0\}$. Obviously, the standard numerical/auto differentiation based on $u_x+u_t$ are computed using values on both sides of the interface $\Sigma$.

%For example, consider the simplest model problem in $\Omega=(0,1)\times (0,1)\subset\R^2$ with $\bm\beta=(1,1)^T$. Let $\bx=(x,t)$ be any point on the discontinuity interface $\Sigma=\{(x,t)\in\Omega : x+t=0\}$. Then the standard numerical/auto differentiation based on $u_x+u_t$ are computed using values on both sides of the interface $\Sigma$.}
%For small $\tau$, the discrete directional differential operator $D_{\scriptstyle{\bm\beta},\tau}$ defined in \cref{finite_diff} ensures that the derivative is computed without crossing the discontinuity interface.}

To define a discrete divergence operator based on \cref{div}, we associate each integration point $\bz$ with a corresponding control volume $K_{\bz}$ that contains the point. The discrete divergence of the total flux $\bF(v)$ at $\bz$ is then given by
\begin{equation}\label{div-dis}
    \div\!\!_{_\cT}\bF(v(\bz)) = \dfrac{1}{\big| K_\bz\big|} \mathcal{Q}_{\partial K_\bz} \left(\bF(v)\!\cdot\! \bn\right),
\end{equation}
where $\mathcal{Q}_{\partial K_\bz}(\cdot)$ denotes a {\it composite quadrature rule} applied over the boundary $\partial K_\bz$ of the control volume $K_\bz$ and $\bn$ is the unit outward vector normal to the boundary $\partial K_\bz$.  

Under the midpoint quadrature rule $\cQ_{K}$, each subdomain $K\in \cT$ has a single integration point $\bz_K=(\bx_K,t_K)$, typically chosen as the centroid of $K$. In this case, the control volume is simple $K_{\bz_K}=K$. 

More generally, suppose the quadrature rule $\cQ_{K}$ for a subdomain $K\in\cT$ has $J$ integration points,
\[
\bz_{K_j}=(\bx_{K_j},t_{K_j})\in K\in\cT, \,\,\mbox{ for }\, j=1,\ldots,J.
\]
Let $\cT_K=\{K_j\}_{j=1}^J$ be a partition of $K$ such that $\bz_{K_j}\in K_j$, where $K_j$ is referred to as the control volume of the integration point $\bz_{K_j}$. Let $\mathcal{Q}_{\partial K_j}(\cdot)$ be a composite quadrature rule over the boundary $\partial K_j$, then the discrete divergence operator $\div\!\!_{_\cT}$ at the integration point $\bz_{K_j}$ can be defined analogously via \cref{div-dis}.

%let $\cT=\{K\}$ be a reasonably fine partition of the computational domain $\Omega\times I$, $\bz_K=(\bx_K,t_K)$ be the centroid of subdomain $K\in \cT$, and $\mathcal{Q}_{\partial K}(\cdot)$ be a {\it composite quadrature rule} over the boundary $\partial K$ of subdomain $K$. Then the discrete divergence operator is defined by\begin{equation}\label{div-dis}    \div\!\!_{_\cT}\bF(u(\bz_K)) = \dfrac{1}{\big| K\big|} \mathcal{Q}_{\partial K} \left(\bF(u)\!\cdot\! \bn\right),\end{equation} where $\bn$ is the unit outward vector normal to the boundary $\partial K$. 

The general definition of the discrete divergence operator $\div\!\!_{_\cT}$ in \cref{div-dis} depends on the quadrature rule $\mathcal{Q}_{\partial K}(\cdot)$ over the boundary $\partial K$, which, in turn, depends on the geometry of $K$. Since the integration mesh $\cT$ is independent of the physical mesh induced by the NN approximation, it is practically advantageous to construct $\cT$ as  a composite mesh generated by the AQR procedure in Algorithm 4.1. In this setup, each element $K\in \cT$ is taken to be simple cell such as a rectangle, a cuboid, or a hypercube in two, three, or higher dimensions, with all faces aligned parallel to the coordinate hyperplanes. For such integration mesh $\cT$ in both two and three dimensions, explicit definitions of $\div\!\!_{_\cT}\bF(u(\bz_K))$ were proposed and rigorously analyzed in \cite{Cai2023nonlinear}, specifically in the context of discontinuous solutions $u$. These formulations ensure that the numerical divergence remains consistent with the underlying physics, even across discontinuities and sharp interface features.

%\subsection{Inflow Boundary and Initial Conditions}\label{s:BIcs}

\subsection{The LSNN Method}\label{s:LSNNs}

Denote the collections of the inflow boundary faces and the initial-time faces of the integration mesh $\cT$ by 
\begin{equation}\label{II}
    \mathcal{E}_{-}=\{E=\partial K\cap \Gamma_{-}:K\in\mathcal{T}\} \,\mbox{ and }\, \mathcal{E}_{0}=\{E=\partial K\cap \left(\Omega\times \{0\}\right):K\in\mathcal{T}\},
    \end{equation} 
respectively. For each face $E$ in $\mathcal{E}_{-}$ or $\mathcal{E}_{0}$, let $\mathcal{Q}_E(w)$ denote a quadrature rule applied to an integrand $w$ defined on $E$. We now define the discrete least-squares functionals. For problem \cref{pde1b}, the functional is defined by 
\begin{equation}\label{dls1}
\hat{\mathcal{L}}_{_{{\cal T}}}\big(v\big)=\!\!\sum_{K\in\mathcal{T}}\!\! \mathcal{Q}_K\!\left((D_{\scriptstyle{\bm{\beta},\tau}}v+\gamma v\!-\!f)^2\right) \,\mbox{ and }\,
\mathcal{L}_{_{{\cal T}}}\big(v;{\bf f}\big)= \hat{\mathcal{L}}_{_{{\cal T}}}\big(v\big) +\!\!\sum_{E\in\mathcal{E}_{-}} \!\! \mathcal{Q}_E\left(|\bm{\beta}\! \cdot \!\bm{n}|(v-g)^2\right)
\end{equation} 
and for problem \cref{pde2b}, it is defined as
\begin{equation}\label{dls2}
\left\{\begin{array}{ll}
\hat{\mathcal{L}}_{_{{\cal T}}}\big(v\big)=\sum\limits_{K\in\mathcal{T}} \mathcal{Q}_K\left(\div\!_{_\cT} \bF(v)\right) \quad\mbox{ and }\, &\\[4mm]
\mathcal{L}_{_{{\cal T}}}\big(v;{\bf f}\big)= \hat{\mathcal{L}}_{_{{\cal T}}}\big(v\big) +\sum\limits_{E\in\mathcal{E}_{-}}\mathcal{Q}_E\left((v-g)^2\right) +\sum\limits_{E\in\mathcal{E}_{0}}\mathcal{Q}_E\left((v-u_0)^2\right) & 
\end{array}\right.
\end{equation} 
To unify the notation across both formulations, let 
\[
\mathcal{V} = \mathcal{V}_{\scriptstyle{\bm\beta}}\,\mbox{ or }\, \cV_{\bff} \quad\mbox{and}\quad \mathcal{V}(h) = \mathcal{V}_{\scriptstyle{\bm\beta}}(g) \,\mbox{ or }\, \cV_{\bff}(\bg)
\]
for problem \cref{pde1b} or \cref{pde2b}, respectively.
The LSNN method then seeks a function $u_{_{N,\!\cT}}\in \cM(l)$ that minimizes the discrete least-squares functional. Specifically, it solves either
% The LSNN method for problem \cref{pde1b} or \cref{pde2b} is to either (1) seek  $u_{_{N,\!\cT}}\in \cM(l)$ such that 
 \begin{equation}\label{discrete_minimization_functional}
  \mathcal{L}_{_{{\cal T}}} \big(u_{_{N,\!\cT}};{\bf f}\big) 
  = \min\limits_{v\in \cM(l)} \mathcal{L}_{_{{\cal T}}}\big(v;{\bf f}\big).
\end{equation} 
or when constrained to the admissible set $\mathcal{V}(h_{_\cT})$, solves,
 % (2) seek $u_{_{N,\!\cT}}\in \cM(l) %\cap \mathcal{V}(h_{_\cT})
% $ such that
\begin{equation}\label{discrete-minimization2b}
    {\mathcal{L}}_{_\cT}\left(u_{_{N,\!\cT}};{\bf f}\right) = \min_{\small v\in \cM(l)\cap \mathcal{V}(h_{_\cT})} {\mathcal{L}}(v;{\bf f}) \quad\mbox{or}\quad \hat{\mathcal{L}}_{_\cT} \left(u_{_{N,\!\cT}}\right) = \min_{\small v\in \cM(l)\cap \mathcal{V}(h_{_\cT})} \hat{\mathcal{L}}(v),
\end{equation} 
where $h_{_\cT}$ is an approximation of the prescribed data $h=g$ or $\bg$.

We now discuss how to weakly enforce the inflow boundary and initial conditions  %The least-squares functionals in \cref{dls1} and \cref{dls2} enforce the inflow boundary and initial conditions through penalization: the summation terms over $E$ in $\mathcal{E}_{-}$ and $\mathcal{E}_{0}$. Below, we impose them 
through the physics-preserved numerical differentiation in \Cref{s:DDO}. 
For simplicity, assume that the quadrature rule $\cQ_K(\cdot)$ is the midpoint rule. In this case, the centroid of $K$, $\bz_K=\bx_K$ or $(\bx_K,t_K)$, is the sole integration point within each subdomain $K$. For each inflow boundary or initial face $E\in \mathcal{E}_{-} \cup \mathcal{E}_{0}$, there exists a subdomain $K\in \cT$ such that $E\in \partial K$. We denote this face by $E_K$ to indicate that $E$ is part of the boundary $\partial K$ of $K$.

For each inflow boundary face $E_K\in \mathcal{E}_{-}$, the directional derivative $D_{\scriptstyle{\bm\beta}, \tau}v(\bx_K)$, as defined in \cref{finite_diff}, is computed by choosing the derivative step size $\tau$ so that the backward point $\bx_K - \tau {\bm{\beta}}(\bx_K)$ lies on the boundary face $\mathcal{E}_{-}$. The directional derivative is then given by 
\begin{equation}\label{finite_diff2}
     D_{\scriptstyle{\bm\beta}, \tau}v(\mathbf{x}_K)= \frac{v(\bx_K)-g\big(\bx_K - \tau {\bm{\beta}}(\bx_K)\big)}{\tau},
\end{equation}
where $g$ is the prescribed inflow boundary condition in \cref{pde1}. Similarly, for problem \cref{pde2}, the discrete divergence operator at the integration point $\bz_K$ is modified to incorporate the boundary or initial data on the relevant face $E_K$. Specifically, the discrete divergence becomes
\begin{equation}\label{div-dis2}
    \div\!\!_{_\cT}\bF(u(\bz_K)) = \left\{\begin{array}{ll}
    \dfrac{1}{\big| K\big|} \left( \mathcal{Q}_{\partial K \setminus E_K} \big(\bF(u)\!\cdot\! \bn\big) + \mathcal{Q}_{E_K} \big(\bF(g)\!\cdot\! \bn\big) \right), & E_K\in \mathcal{E}_{-},\\[4mm]
    \dfrac{1}{\big| K\big|} \left( \mathcal{Q}_{\partial K \setminus E_K} \big(\bF(u)\!\cdot\! \bn\big) + \mathcal{Q}_{E_K} \big(\bF(u_0)\!\cdot\! \bn\big) \right), & E_K\in \mathcal{E}_{0}.
    \end{array}
    \right.
\end{equation}
where $g$ and $u_0$ denote inflow boundary and initial conditions in \cref{pde2}, respectively.

\begin{comment}
Denote the modified least-squares functionals by
\begin{equation}\label{m-lsf}
    \mathcal{G}_{_{{\cal T}}}\big(v;{\bf f}\big)=\sum_{K\in\mathcal{T}}\mathcal{Q}_K\left((D_{\bm{\beta}}v +\gamma v-f)^2\right) \quad\mbox{and}\quad
    \mathcal{G}_{_{{\cal T}}}\big(v;{\bf f}\big)=\sum_{K\in\mathcal{T}}\mathcal{Q}_K\left(\div\!_{_\cT} \bF(v)\right)
\end{equation}
for problems \cref{pde1} and \cref{pde2}, respectively, where the discrete directional and divergence operators at subdomains $K\in\cT$, whose boundary intersects $\mathcal{E}_{-}$ or $\mathcal{E}_{0}$, are modified in the respective \cref{finite_diff2} and \cref{div-dis2}. Then the modified least-squares least-squares (LSNN) method for problems \cref{pde1} or \cref{pde2} is to seek  $u_{_{N,\!\cT}}\in \cM(l)$ such that
\begin{equation}\label{discrete_minimization_functional2}
  \mathcal{G}_{_{{\cal T}}} \big(u_{_{N,\!\cT}};{\bf f}\big) 
  = \min\limits_{v\in \cM(l)} \mathcal{G}_{_{{\cal T}}}\big(v;{\bf f}\big).
\end{equation}
\end{comment}

\section{Efficient and Reliable Iterative Solvers}\label{s:solver}

Both $\mathcal{L}_{_{{\cal T}}}\big(v;{\bf f}\big)$ and $\hat{\mathcal{L}}_{_{{\cal T}}}\big(v\big)$ are convex functionals with respect to $v$, but become non-convex when viewed as functions of the NN parameters. As a result, the discrete problems formulated in \cref{discrete_minimization_functional} or \cref{discrete-minimization2b} are non-convex optimization problems over the NN parameters. These high-dimensional, non-convex optimization problems are typically computationally intensive and complex, and they remain a major bottleneck in applying NNs to numerically solve PDEs. Despite these challenges, there has been active research and some notable progress in developing more efficient and reliable iterative solvers (training algorithms) and in designing effective initialization strategies for NN-based PDE methods \cite{SgGN, CaiDokFalHer2024a, CaiDokFalHer2024b}. %This section discusses some algebraic structures of the discretization that may be used for designing efficient and reliable iterative solvers

The Adam (Adaptive Moment Estimation) method \cite{kingma2015} is a first-order gradient-based optimization algorithm widely used in training NNs. At each iteration, Adam updates the NN parameters using an adaptive learning rate based on estimates of the first and second moments of the gradients. This adaptive learning rate mechanism makes it particularly effective for large-scale optimization problems with noisy or sparse gradients. However, for solving high-dimensional, non-convex optimization problems arising from NN-based PDE discretizations, first-order methods like Adam often suffer from slow convergence, motivating the development of more efficient second-order solvers that can better exploit the structure of the optimization landscape.

Below, we first explore the algebraic structure of the problem and then introduce a second-order Gauss-Newton based iterative method. As a nonlinear PDE, \cref{pde2b} possesses its own intrinsic nonlinearity that warrants special treatment. In this section, we focus solely on the linear problem in \cref{pde1b}. For simplicity of presentation, we use the second minimization problem defined in \cref{discrete-minimization2b}, namely the task of finding $u^*_{_{N,\!\cT}} (\bx)\in \cM(l)$ such that
\begin{equation}\label{d-min-hat}
    u^*_{_{N,\!\cT}}\left(\bx\right) = \argmin_{v\in \cM(l)\cap \mathcal{V}(h_{_\cT})} \hat{\mathcal{L}}(v)
\end{equation}
where $\hat{\mathcal{L}}(v)$ defined in \cref{dls1}, as an example to illustrate the algebraic structures underlying  this non-convex optimization problem. These structures can be leveraged for designing more efficient and reliable iterative solvers.

Problem \cref{d-min-hat} is often called separable nonlinear least-squares (SNLS) problem (see, e.g., \cite{Kaufman75}), due to the presence of both the linear and nonlinear parameters. There are two approaches for solving SNLS problems: (1) block iterative methods that alternate between updates of the linear and the nonlinear parameters (serving as outer iteration), and (2) the Variable Projection (VarPro) method of Golub-Pereyra \cite{Golub1973}, which eliminates the linear parameters to reduce the dimensionality of the problem. Although the VarPro method yields a problem with fewer degrees of freedom, it alters the original nonlinear structure of the SNLS problem. Therefore, in this work, we will focus exclusively on the first approach.

To this end, any function $v(\bx)\in \cM(l)$ has the following form 
\begin{equation}\label{uNT}
v(\bx) = %u_{_{N,\!\cT}}(\bx; \hat{\bc}, {\small \bTheta})= 
c_0+\sum_{i=1}^{n_{l}} c_i \sigma\left(\bomega^{(l)}_i\bx^{(l-1)}+b^{(l)}_i\right),
\end{equation}
where $\bx^{(l-1)}$ is recursively defined in \cref{layerdef}. Clearly, $v(\bx)$ is determined by the linear and nonlinear parameters
\[
\hat{\bc} =\left(c_0,{\bc}\right)= \left(c_0,c_1,\ldots,c_{n_l}\right) \in \R^{n_l+1} \,\mbox{ and }\,  {\small \bTheta} = \left\{ \br^{(k)} \right\}_{k=1}^l = \left\{ \left(\br^{(k)}_1,\ldots, \br^{(k)}_{n_k} \right)^T\right\}_{k=1}^l,
\]
respectively, where $\br^{(k)}_i=\left(b^{(k)}_i, \bomega^{(k)}_i \right)$ is the bias and weights of the $i^{th}$ neuron at the $k^{th}$ hidden layer. To indicate this  dependence, we often write $v(\bx)= v(\bx; \hat{\bc}, {\small \bTheta})$. Now, a solution 
\[
u^*_{_{N,\!\cT}}(\bx)=u_{_{N,\!\cT}}(\bx; \hat{\bc}^*, {\small \bTheta}^*)
\]
of the problem in $\cref{d-min-hat}$ is to minimize the least-squares functional $\hat{\mathcal{L}}(v)$ over $\cM(l)\cap \mathcal{V}(h_{_\cT})$; or equivalently, the optimal parameters $(\hat{\bc}^*, {\small \bTheta}^*)$ minimize 
the function $\hat{\mathcal{L}}\left(v(\cdot; \hat{\bc}, {\small \bTheta})\right)$ of variables $(\hat{\bc}, {\small \bTheta})$ in a domain of high dimension. Hence, $(\hat{\bc}^*, {\small \bTheta}^*)$ satisfies the following optimality conditions
\begin{equation}\label{OC}
    \nabla_{\bc} \hat{\mathcal{L}}_{_{{\cal T}}} \left(u_{_{N,\!\cT}}(\cdot; \hat{\bc}^*, {\small \bTheta}^*)\right)  ={\bf 0}
    \quad\mbox{and}\quad
    \nabla_{\footnotesize \bTheta} \hat{\mathcal{L}}_{_{{\cal T}}} \left(u_{_{N,\!\cT}}(\cdot; \hat{\bc}^*, {\small \bTheta}^*)\right)  ={\bf 0},
\end{equation}
where $\nabla_{\hat{\bc}}$ and $\nabla_{\footnotesize \bTheta}$ denote the gradients with respect to $\hat{\bc}$ and ${\small\bTheta}$, respectively.

Below we discuss algebraic structures of \cref{OC}. First, this is a coupled nonlinear system of algebraic equations on $\hat{\bc}^*$ and ${\small \bTheta}^*$. Next, we derive specific forms of these algebraic equations. To do so, let $\sigma_0(\bx)=1$ and $\sigma_i(\bx)=\sigma\left(\bomega^{(l)}_i\bx^{(l-1)}+b^{(l)}_i\right)$ for $i=1,\ldots,n_l$, and set 
\[
\bSigma(\bx)=(\sigma_1(\bx),\ldots,\sigma_{n_l}(\bx))^T \quad\mbox{and}\quad \hat{\bSigma}(\bx)=(\sigma_0(\bx),\sigma_1(\bx),\ldots,\sigma_{n_l}(\bx))^T.
\]
Since the ReLU activation function $\sigma(s)$ is not point-wisely differentiable, it is then difficult to compute gradients of the discrete LS functional $\hat{\mathcal{L}}_{_{{\cal T}}} \left(v(\cdot; \hat{\bc}, {\small \bTheta})\right)$. Instead, we introduce an intermediate functional before numerical integration, \begin{equation}\label{inter-L-hat}
    \hat{\mathcal{L}}_{\tau} \left(v(\cdot; \hat{\bc}, {\small \bTheta})\right)= \int_\Omega (D_{\scriptstyle{\bm{\beta},\tau}}v+\gamma v\!-\!f)^2\,d\bx,
\end{equation}
compute gradients of this functional, and then approximate the coefficient matrix and the right-hand side vector by numerical quadrature.

Because $\sigma(s)$ has the first order weak derivative, by the facts that $v(\bx)= \hat{\bSigma}(\bx)^T \hat{\bc}$ and that $\nabla_{\hat{\bc}} v(\bx; \hat{\bc}, {\small \bTheta}) = \hat{\bSigma}(\bx)$, we have
\begin{eqnarray*}
&& \dfrac12 \nabla_{\hat{\bc}} \hat{\mathcal{L}}_\tau \left(v(\cdot; \hat{\bc}, {\small \bTheta})\right) = \int_\Omega \!\left(D_{\scriptstyle{\bm\beta},\tau} \,v +\gamma\, v-f\right) \nabla_{\hat{\bc}} \left(D_{\scriptstyle{\bm\beta},\tau} +\gamma\right)v\, d\bx \\[2mm] & = &\!\!\!\! \int_\Omega  \!\left(D_{\scriptstyle{\bm\beta},\tau} \,v +\gamma\, v-f\right)  \left(D_{\scriptstyle{\bm\beta},\tau} +\gamma\right) \hat{\bSigma}(\bx) d\bx = \!\int_\Omega \!\left(D_{\scriptstyle{\bm\beta},\tau} +\gamma\right) \hat{\bSigma}(\bx) \left\{\left(D_{\scriptstyle{\bm\beta},\tau} +\gamma\right) \hat{\bSigma}(\bx)^T \hat{\bc} -f \right\}  d\bx,
\end{eqnarray*}
which, together with the optimality condition, implies
\begin{equation}\label{OC1}
   {\bf 0} = \dfrac12 \nabla_{\hat{\bc}} \hat{\mathcal{L}}_\tau \left(v(\cdot; \hat{\bc}, {\small \bTheta})\right) =  \int_\Omega \left\{\left[(D_{\scriptstyle{\bm\beta},\tau} +\gamma) \hat{\bSigma}\right]\left[(D_{\scriptstyle{\bm\beta},\tau} +\gamma) \hat{\bSigma}^T\right] \hat{\bc} -f  \left[(D_{\scriptstyle{\bm\beta},\tau} +\gamma) \hat{\bSigma}\right]\right\}  d\bx.
\end{equation}

Approximating the integral by numerical quadrature, then the first equation in \cref{OC} implies the following system of algebraic equations 
\begin{equation}\label{linear system}
    \bm{A}\left({\small \bTheta}\right)\, \hat{\bc} = F\left({\small \bTheta}\right),
\end{equation}
where $\bm{A}\left({\small \bTheta}\right)$ and $F\left({\small \bTheta}\right)$ are the coefficient matrix of order $(n_l+1)\times (n_l+1)$ and the right-hand side vector $(n_l+1)\times 1$ given by 
\begin{equation}\label{A-F}
\left\{\begin{array}{l}
    \bm{A}\left({\small \bTheta}\right)= \sum\limits_{K\in\mathcal{T}}\mathcal{Q}_K\left(\left[(D_{\scriptstyle{\bm\beta},\tau} +\gamma) \hat{\bSigma}\right]\left[(D_{\scriptstyle{\bm\beta},\tau} +\gamma) \hat{\bSigma}^T\right]\right)  \quad \mbox{and}\quad\\ [6mm]
 F\left({\small \bTheta}\right)=\sum\limits_{K\in\mathcal{T}}\mathcal{Q}_K\left(f \left[(D_{\scriptstyle{\bm\beta},\tau} +\gamma) \hat{\bSigma}\right] \right),
\end{array}\right.
\end{equation}
respectively. Here the actions of the numerical integration and numerical differentiation operators $\mathcal{Q}_K$ and $D_{\bm{\beta},\tau}$ are applied component-wisely. 

At each step of a block iterative method for solving the SNLS problem in \cref{OC}, given the current estimate of the nonlinear parameters ${\small \bTheta}$, \cref{linear system} is a linear system with the coefficient matrix $\bm{A}\left({\small \bTheta}\right)$. Let  $a_{ij}\left({\small \bTheta}\right)$ be the $ij$-element of $\bm{A}\left({\small \bTheta}\right)$, then 
\[
a_{ij}\left({\small \bTheta}\right) = \sum\limits_{K\in\mathcal{T}} \mathcal{Q}_K\left(\left[D_{\bm{\beta},\tau}\sigma_i +\gamma \sigma_i\right] \left[D_{\bm{\beta},\tau}\sigma_j +\gamma \sigma_j\right]\right) =a_{ji}\left({\small \bTheta}\right),
\]
which shows that $\bm{A}\left({\small \bTheta}\right)$ is symmetric. Under reasonable assumptions (e.g., the assumption for one hidden layer is that the corresponding breaking hyperplanes are distinct), 
the set $\left\{ \sigma_i (\bx)\right\}_{i=1}^{n_l}$ is linearly independent, implying that $\bm{A}\left({\small \bTheta}\right)$ is positive definite (see, e.g., \cite{SgGN} for shallow NN). However, due to the global support of the basis functions $\left\{ \sigma_i (\bx)\right\}_{i=1}^{n_l}$, $\bm{A}\left({\small \bTheta}\right)$ is a dense matrix and can be highly ill-conditioned (see, e.g., \cite{CaiDokFalHer2024a, CaiDokFalHer2024b}, especially for shallow NNs in one dimension). Consequently, gradient-based optimization methods may suffer from slow convergence, limiting their efficiency for solving such systems.
%A large condition number implies that the classic iterative methods such as the Richardson method (equivalently, the gradient descent methods) reduces the algebraic error super slowly. 

Similar systems of algebraic equations to \cref{A-F} arise from the least-squares approximation using shallow ReLU NNs \cite{SgGN} and the shallow Ritz method for one-dimensional diffusion and diffusion-reaction problems \cite{CaiDokFalHer2024a, CaiDokFalHer2024b}. In those special cases, efficient and reliable iterative solvers have been developed and analyzed. However, in broader NN applications, the design of fast solvers for the linear parameters remains a critical and largely unresolved challenge. When the number of the linear parameters is relatively small, techniques such as truncated singular value decomposition (SVD) can effectively mitigate the issues caused by large condition numbers and provide accurate solutions despite the ill-conditioning \cite{SgGN}. Nevertheless, scalable and efficient algorithms for larger systems are still an active area of research. 

\begin{comment}
\begin{lemma}
    The coefficient matrix $\bm{M}\left({\small \bTheta}^{(l)}\right)$ defined in \cref{M-F} is symmetric. For the fixed nonlinear parameters ${\small \bTheta}^{(l)}$, assume that functions $\left\{\varphi_{i}^{(l)}(\bx)\right\}_{i=0}^{n_l}$ are linearly independent and that a discrete version of \cref{gamma} holds: there exists a positive constant $\hat{\gamma}_0$ such that
    , then $\bm{M}\left({\small \bTheta}^{(l)}\right)$ is positive definite.
\end{lemma}

\begin{proof}
The symmetry of $\bm{M}\left({\small \bTheta}^{(l)}\right)$ is obvious. For any $\bxi\in \R^{n_l+1}$, let $v(\bx)= \bxi^T\bvarphi(\bx) \in \cM(l)$, then we have 
\[
\bxi^T \bm{M}\left({\small \bTheta}^{(l)}\right) \bxi =\sum\limits_{K\in\mathcal{T}}\mathcal{Q}_K\left(\left(\bxi^T\left[D_{\bm{\beta}}\bvarphi^{(l)} +\gamma \bvarphi^{(l)}\right] \right)^2\right) = \sum\limits_{K\in\mathcal{T}}\mathcal{Q}_K\left(\left[D_{\bm{\beta}}v+\gamma v\right]^2\right)
\]
which, together with the assumption, implies positive definiteness of $\bm{M}^{(l)}\left({\small \bTheta}^{(l)}\right)$.
\end{proof}
\end{comment}

At each step of the block iterative method, the second equation in \cref{OC} remains a nonlinear system of algebraic equations for ${\small \bTheta}$, given the current estimate of the linear parameters $\hat{\bc}$. A natural choice for solving such NLS problems is the Gauss-Newton (GN) method \cite{dennis1996numerical,ortega2000iterative}, since its associated GN matrix is always positive semi-definite. A major limitation of the GN method is the potential singularity of the GN matrix, which often necessitates regularization strategies, most notably the shifting technique used in the Levenberg-Marquardt (LM) method \cite{Levenberg, Marquardt}, to ensure invertibility. 
These regularization techniques, while widely adopted, modify the original optimization problem and introduce additional complexity. In particular, selecting an appropriate shift parameter is highly problem-dependent and remains a challenging and delicate task in practice.

For the least-squares approximation of a target function using shallow NNs, in \cite{SgGN} we addressed the limitations of the LM by deriving a structured form of the GN matrix and explicitly eliminating its singularity. The resulting structure-guided Gauss-Newton (SgGN) method achieves faster convergence and significantly greater accuracy than the LM approach. An extension of the SgGN method to least-squares approximation problems using deep ReLU NNs will be presented in a forthcoming paper.

%one may employ the commonly used first-order gradient-based methods (see, e.g., survey papers \cite{SIAM_Review18, Survey19, SCZZ2020_Survey}), second-order methods (see, e.g., survey papers \cite{SIAM_Review18, Survey19, SCZZ2020_Survey}), or the Gauss-Newton (GN) method \cite{dennis1996numerical,ortega2000iterative} for nonlinear least-squares optimization. 
\begin{comment}

For the shallow ReLU NN with $n$ neurons, the nonlinear parameters consist of 
\[
%\left\{ \left(b_i, \bomega_i \right) \right\}_{i=1}^{n} = \left\{ \br_i\right\}_{i=1}^{n} 
{\small \bTheta}=\br = (\br_1,\ldots, \br_n)^T \quad\mbox{with }\, \br_i= \left(b_i, \bomega_i \right),
\]
where the superscript indicating the hidden-layer number is omitted when there is only one hidden-layer. A solution of problem \cref{discrete_minimization_functional} or \cref{discrete_minimization_functional2} has the form of
\[
u_{_{N,\!\cT}}(\bx) = \sum_{i=0}^{n} c_i \sigma_{i}(\bx) =\bc^T \bSigma(\bx;\br),
\]
where $\bSigma(\bx;\br)= \left(\sigma_0(\bx),\ldots, \sigma_n(\bx;\br)\right)^T$ and 
\[
\sigma_{0}(\bx)=1
\quad\mbox{and}\quad 
\sigma_{i}(\bx;\br)= \sigma  \left(b_i +\bomega_i\cdot\bx \right) = \sigma  \left(\br_i\cdot\by \right) \mbox{ with } \by^T= (1,\bx) .
\]
The system of linear equations in \cref{linear system} becomes
\begin{equation}\label{linear system2}
    \bm{A}(\br) \bc = F(\br),
\end{equation}
where the coefficient matrix $\bm{A}(\br)$ and the right-hand side $F(\br)$ are defined accordingly as in \cref{A-F}.
\end{comment}

Next, we describe the SgGN method for solving the second equation in \cref{OC} when using shallow ReLU NNs. Any NN function $v\left(\bx\right)\in \cM(1)$ has the form of
\begin{equation}\label{v1a}
v\left(\bx\right) = v(\bx;\hat{\bc},\br) = c_0+\sum\limits_{i=1}^{n}c_i \sigma\left(\bomega_i\bx+b_i\right).
\end{equation}
Here and thereafter, we drop the subscript $1$ of $n_1$ and the superscript $(1)$ for notational simplicity. The nonlinear parameters consist of $\br = (\br_1,\ldots, \br_n)$, with each $\br_i=\left(b_i, \bomega_i \right)$.  Substituting this into the second equation in \cref{OC}, we obtain
\begin{equation}\label{OC2}
    {\bf 0} = \nabla_\br \hat{\mathcal{L}}_{_{{\cal T}}} \left(u_{_{N,\!\cT}}(\cdot; \hat{\bc}^*, \br^*)\right)
\end{equation}
Again, we use the intermediate functional $\hat{\mathcal{L}}_\tau$ to derive a structured form of the GN matrix and then approximate the integral by quadrature.

Denote by $H(s)=\sigma^\prime(s)=\left\{\begin{array}{cc}
    0, & s<0, \\
    1, & s>0
\end{array}\right.$ %for $s<0$ and $H(s)=\sigma^\prime(s)=1$ for $s>0$ 
the Heaviside step function, where $\sigma^\prime(s)$ is the weak derivative of $\sigma(s)$ in the distribution sense. Let $I_{d+1}$ be the order-$(d+1)$ identity matrix, $D(\bc)$ be the diagonal matrix with the $i^{\mbox{\scriptsize th}}$-diagonal elements $c_i$, and $\by=(1,x_1,\ldots,x_d)^T$ be the homogeneous-coordinates. A direct calculation (see (4.5) in \cite{SgGN}) gives
\begin{equation}\label{grad-r}
   \nabla_{\br} v(\bx;\hat{\bc},\br) 
   = \left(D(\bc)\otimes I_{d+1}\right) \left(\bH(\bx)\otimes \by\right),
\end{equation}
where $\bH(\bx)=\big(H_1\left(\bx;\br_1\right), \ldots, H_n\left(\bx;\br_n\right)\big)^T$ with $H_i\left(\bx;\br_i\right)= H\left(\bomega_i\bx+b_i\right) = H\left(\br_i\by\right)$ and the symbol $\otimes$ denotes the Kronecker product of two matrices/vectors. Denote by
\begin{equation}\label{G}
    \bG(\hat{\bc},\br)=\int_\Omega \left(D_{\scriptstyle{\bm\beta},\tau} \,v +\gamma\, v-f\right)  \left(D_{\scriptstyle{\bm\beta},\tau} +\gamma\right)  \left(\bH(\bx)\otimes \by\right) \,d\bx
\end{equation}
a scaled gradient vector of $\hat{\mathcal{L}}_\tau \left(v(\cdot; \hat{\bc}, \br)\right)$ with respect to $\br$.
By \cref{grad-r}, we have
\[
\nabla_{\br} \hat{\mathcal{L}}_\tau \left(v(\cdot; \hat{\bc}, \br)\right) = 2\int_\Omega \!\left(D_{\scriptstyle{\bm\beta},\tau} \,v +\gamma\, v-f\right)  \left(D_{\scriptstyle{\bm\beta},\tau} +\gamma\right) \nabla_{\br} v\, d\bx = 2\left(D(\bc)\otimes I_{d+1}\right) \bG(\hat{\bc},\br) 
\]
The principal part of the Hessian without the second-order derivative with respect to $\br$ is called the GN matrix given by
\[
\left\{\begin{array}{l}
\mathcal{G}(\bc,\br) = 2 \left(D(\bc)\otimes I_{d+1}\right)  \mathcal{H}(\br) \left(D(\bc)\otimes I_{d+1}\right)\\[2mm]
\mbox{with }\, \mathcal{H}(\br) = \int_\Omega \left[ \left(D_{\scriptstyle{\bm\beta},\tau} +\gamma\right)  \left(\bH(\bx)\otimes \by\right)\right] \left[ \left(D_{\scriptstyle{\bm\beta},\tau} +\gamma\right)  \left(\bH(\bx)\otimes \by\right)^T\right] d\bx,
\end{array}\right.
\]
where $\mathcal{H}(\br)$ is referred to as the layer GN matrix.
Approximating the integral by numerical quadrature, we then have the corresponding components 
\begin{equation}\label{NL-T}
\left\{\begin{array}{l}
 \bG_{_{{\cal T}}}(\hat{\bc},\br)=\sum\limits_{K\in \cT} \mathcal{Q}_K\left\{\left(D_{\scriptstyle{\bm\beta},\tau} \,v +\gamma\, v-f\right)  \left(D_{\scriptstyle{\bm\beta},\tau} +\gamma\right)  \left(\bH\otimes \by\right)\right\}, \\[4mm]
 %\nabla_{\br} \hat{\mathcal{L}}_{_{{\cal T}}} \left(v(\cdot; \hat{\bc}, \br)\right) = 2\left(D(\bc)\otimes I_{d+1}\right) \bG_{_{{\cal T}}}(\hat{\bc},\br),  \\[3mm]
 \mathcal{H}_{_{{\cal T}}}(\br) = \sum\limits_{K\in \cT} \mathcal{Q}_K\left\{ \left[ \left(D_{\scriptstyle{\bm\beta},\tau} +\gamma\right)  \left(\bH\otimes \by\right)\right] \left[ \left(D_{\scriptstyle{\bm\beta},\tau} +\gamma\right)  \left(\bH\otimes \by\right)^T\right]\right\}, \,\mbox{ and }\, 
 \\[4mm]
 \mathcal{G}_{_{{\cal T}}}(\bc,\br) = 2 \left(D(\bc)\otimes I_{d+1}\right)  \mathcal{H}_{_{{\cal T}}}(\br) \left(D(\bc)\otimes I_{d+1}\right).
\end{array}\right.
\end{equation}

Now, we describe the SgGN method \cite{SgGN} for solving problem in \cref{d-min-hat}. The algorithm starts with an initial function approximation $u_{n,_\cT}^{(0)}(\bx)=u_{n,_\cT}\left(\bx;\hat{\bc}^{(0)}, \br^{(0)}\right)$ by  initializing the nonlinear parameters $\br^{(0)}$, as they define the breaking hyperplanes that partition the domain, and effectively forming a computational ``mesh'' for our approximation. Given $\br^{(0)}$, we compute the optimal linear parameters $\hat{\bc}^{(0)}$ on the current physical partition by solving 
\begin{equation}\label{initial-c}
 \bm{A}\left(\br^{(0)}\right)\, \hat{\bc}^{(0)} = F\left(\br^{(0)}\right),
\end{equation} 
where $\bm{A}\left(\br^{(0)}\right)$ and $F\left(\br^{(0)}\right)$ are defined in \cref{A-F}.
For a detailed discussion of this initialization strategy, see \cite{LiuCai1, Cai2021linear}.
Then, given the approximation $u_{n,_\cT}^{(k)}(\bx)=u_{n,_\cT}\left(\bx;\hat{\bc}^{(k)}, \br^{(k)}\right)$ at the $k^{th}$ iteration, the process of obtaining 
\[
u_{n,_\cT}^{(k+1)}(\bx)=u_{n,_\cT}\left(\bx;\hat{\bc}^{(k+1)}, \br^{(k+1)}\right)
\]
proceeds as follows:
\begin{itemize}
\item[(i)] First, identify a set of {\em active} neurons in the current approximation $u_{n,_\cT}^{(k)}(\bx)$ by defining
\begin{equation}
\mathcal{I}_{\text{active}} = \left\{ i \in \{1, \dots, n\} \left|\; \left|{c}_i^{(k)}\right| \right.\ge \epsilon_{\bc} \right\},\label{eq:active}
\end{equation}
where $\epsilon_{\bc}$ is a prescribed tolerance.
The parameters corresponding to these active neurons, denoted as $\tilde{\bc}^{(k)}$ and $\tilde{\br}^{(k)}$, are extracted to form a reduced system. %for refinement.
\item[(ii)] Next, %construct and solve a reduced GN system using only the parameters associated with the active set identified in step (i). Unlike traditional methods that artificially modify the entire system, 
we form $\tilde{\bG}_{_{{\cal T}}}\left(\hat{\bc}^{(k)},\tilde{\br}^{(k)}\right)$, $D\left(\tilde{\bc}^{(k)}\right)$, and $\tilde{\mathcal{H}}_{_{{\cal T}}}\left(\tilde{\br}^{(k)}\right)$, compute the search direction in the reduced space by
\begin{equation}\label{GN2}
{\tilde{\bp}}^{(k+1)} = \left({D}^{-1} \left({\tilde{\bc}}^{(k)}\right)\otimes I_{d+1}\right)  {\tilde{\mathcal{H}}}^{-1}\left({\tilde{\br}}^{(k)}\right)\tilde{\bG}_{_{{\cal T}}}\left(\hat{\bc}^{(k)},\tilde{\br}^{(k)}\right),
\end{equation}
and then map it back to the full parameter space by initializing $\bp^{(k+1)} = \bm{0}$ and setting $\bp^{(k+1)}_i = \tilde{\bp}^{(k+1)}_i$ only for indices $i \in \mathcal{I}_{\text{active}}$. %This selective update strategy explicitly excludes inactive parameters, fundamentally addressing the mathematical structure of the problem rather than using shifting.
\item[(iii)] Then, the nonlinear parameter is updated by
\begin{equation*}
{\br}^{(k+1)} = {\br}^{(k)} - \gamma_{k+1} {\bp}^{(k+1)},
\end{equation*}
where the optimal step size $\gamma_{k+1}$ is computed by minimizing one dimensional function:
\begin{equation*}
\gamma_{k+1} = \argmin_{\gamma \in {\R^{+}_{0}} } \,\hat{\mathcal{L}}_{_\cT}\left(u_{n,_\cT}\left(\cdot;{\hat{\bc}}^{(k)},{\br}^{(k)}-\gamma {\bp}^{(k+1)}\right)\right).
\end{equation*}
%which leads to the n\textbf{}onlinear parameter update\begin{equation*}{\br}^{(k+1)} = {\br}^{(k)} - \gamma_{k+1} {\bp}^{(k+1)}.\end{equation*}
\item[(iv)] Finally, $u_{n,_\cT}^{(k+1)}(\bx)=u_{n,_\cT}\left(\bx;\hat{\bc}^{(k+1)},\br^{(k+1)}\right)$ is obtained by solving
\begin{equation}\label{GN1}
    {\bm A}\left(\br^{(k+1)}\right) \,\hat{\bc}^{(k+1)}=F\left(\br^{(k+1)}\right).
\end{equation} 
%These updated parameters $\hat{\bc}^{(k+1)}$ and ${\br}^{(k+1)}$ together define the improved function approximation $u_{n,_\cT}^{(k+1)}(\bx)=u_{n,_\cT}\left(\bx;\hat{\bc}^{(k+1)},\br^{(k+1)}\right)$.
\end{itemize}

We conclude this section with a few remarks on initialization. The optimization problem in \cref{d-min-hat} is inherently non-convex, making initialization a critical factor for the success of any optimization/iterative/training scheme. This challenge can be addressed by leveraging (1) the physical interpretations of the linear and nonlinear parameters and (2) method of various continuations. 

%The resulting discrete problems in \cref{discrete_minimization_functional} and \cref{discrete_minimization_functional2} are non-convex optimization, and hence initialization is critical for the success of any optimization/iterative/training scheme. The initialization issue may be addressed through (1) the physical meaning of the linear and nonlinear parameters and (2) method of various continuations.  

For the shallow ReLU neural network, since the breaking hyper-planes of neurons form a partition of the computational domain, initialization of the nonlinear parameters $\br$ is given by lying those hyper-planes that uniformly partition the domain. Initialization of the linear parameters $\bc$ is then the solution of \cref{initial-c} %with fixed $\br$ that is a linear problem 
(see \cite{Cai2021linear, Cai2023nonlinear, SgGN, CaiDokFalHer2024a, CaiDokFalHer2024b}).

The adaptive neuron enhancement (ANE) method introduced in \cite{LiuCai1, Cai2023AI} provides a natural method of continuation. The method of model continuation for linear advection-reaction problems with variable advection field was studied in \cite{Cai2021linear}. Finally, the method of subdomain continuation for the block space-time LSNN method was introduced in \cite{Cai2023nonlinear} for the nonlinear hyperbolic conservation laws. 

%To narrow down domain of the nonlinear parameters $\br_i=\left(b_i, \bomega_i \right)$ for possible fewer local/global minimizers, we normalize the weights of each neuron (see, e.g., \cite{Cai2021linear, LiuCai1, SgGN}) by setting $\bomega_i\in \mathcal{S}^{d-1}=\left\{\bxi\in \R^d :\, |\bxi|=1  \right\}$ (the unit sphere in $\R^d)$, where $|\bxi|$ denotes the magnitude of $\bxi$. 

\section{Numerical Experiment}\label{s:NE}
In this section, we present three numerical examples to demonstrate the performance of the LSNN method for linear and nonlinear hyperbolic problems. In each experiment, the discrete LS functionals were minimized using the Adam first-order optimization algorithm \cite{kingma2015}. Implementation of the second-order Gauss-Newton method, as presented in Section 5, is beyond the scope of the current experiments and is deferred to future work or for readers to explore. The structure of the ReLU NN used is denoted as $d$-$n_1$-$n_2\cdots n_{l-1}$-$d_{o}$ for a $l$-layer network, where $n_1$, $n_2$ and $n_{l-1}$ represent the number of neurons in the first, second, and $(l-1)$th layers, respectively. Here, $d$ and $d_{o}$ indicate the input and output dimensions of the problem.

\subsection{A 2D linear problem with a variable advection velocity field} \label{sec6.1}
Consider a variable advective velocity field $\bm{\beta}(x,y) = (1,2x),\,\, (x,y)\in\Omega =(0,1)\times (0,1)$, and the boundary of the input of the problem is $\Gamma_{-}=\{(0,y):y\in(0,1)\}\cup\{(x,0):x\in(0,1)\}$. The inflow boundary condition is given by
\begin{equation*}
g(x,y)=\left\{ \begin{array}{rl}
 y+2,& (x,y)\in \Gamma^1_-\equiv\{(0,y): y\in[\frac{1}{5},1)\}, \\[2mm]
 (y-x^2)e^{-x}, & (x,y)\in \Gamma^2_-=\Gamma_-\setminus \Gamma_-^1.
\end{array}\right.
\end{equation*} 
The exact solution of this linear advection-reaction problem is
\begin{equation}
u(x,y)=\left\{ \begin{array}{rl}
 (y-x^2)e^{-x},& (x,y)\in \Omega_1\equiv\{(x,y)\in\Omega:y< x^2+\frac{1}{5}\}, \\[2mm]
 (y-x^2+2)e^{-x}, & (x,y)\in \Omega_2=\Omega\setminus\Omega_1.
\end{array}\right.
\end{equation}

\begin{figure}[htbp]\label{test5 figure}
\centering
\subfigure[The interface\label{interface5}]{
\begin{minipage}[t]{0.4\linewidth}
\centering
\includegraphics[width=2in]{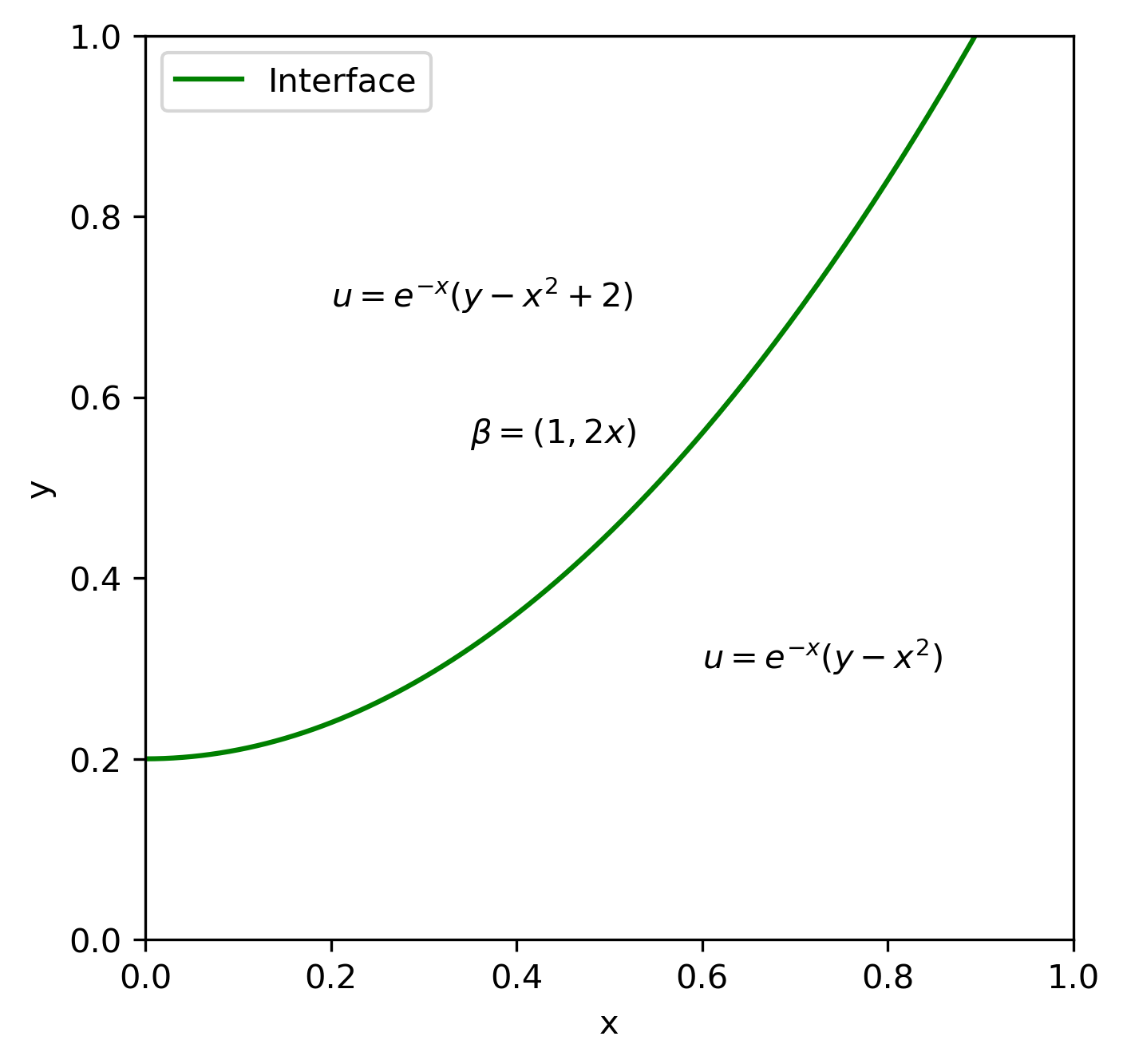}
\end{minipage}%
}%
\hspace{0.2in}
\subfigure[The trace of Figure \ref{comparison5} on $y=1-x$\label{vertical5}]{
\begin{minipage}[t]{0.4\linewidth}
\centering
\includegraphics[width=1.9in]{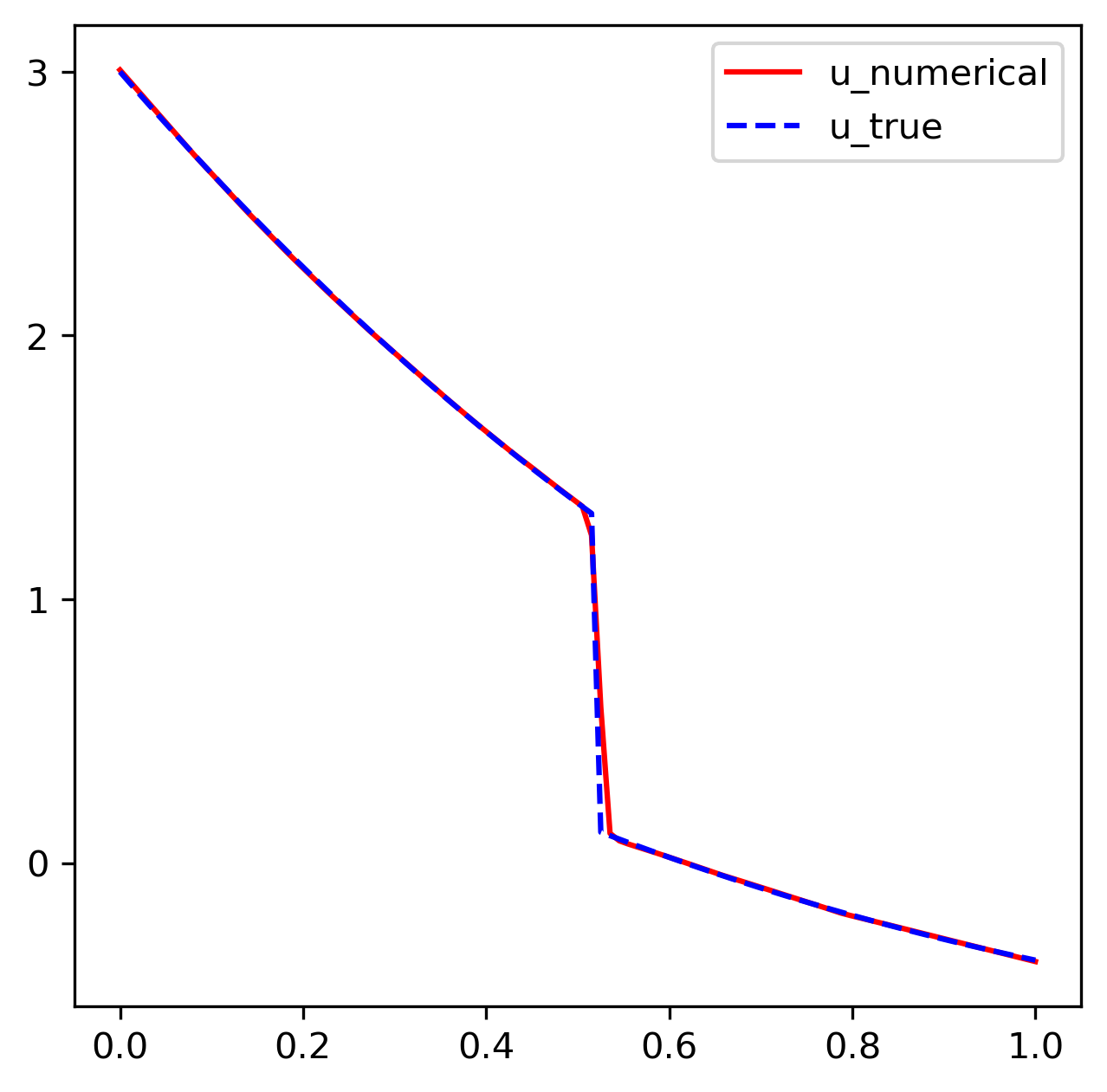}
\end{minipage}%
}%
\\
\subfigure[The exact solution\label{comparison_exact5}]{
\begin{minipage}[t]{0.4\linewidth}
\centering
\includegraphics[width=2.1in]{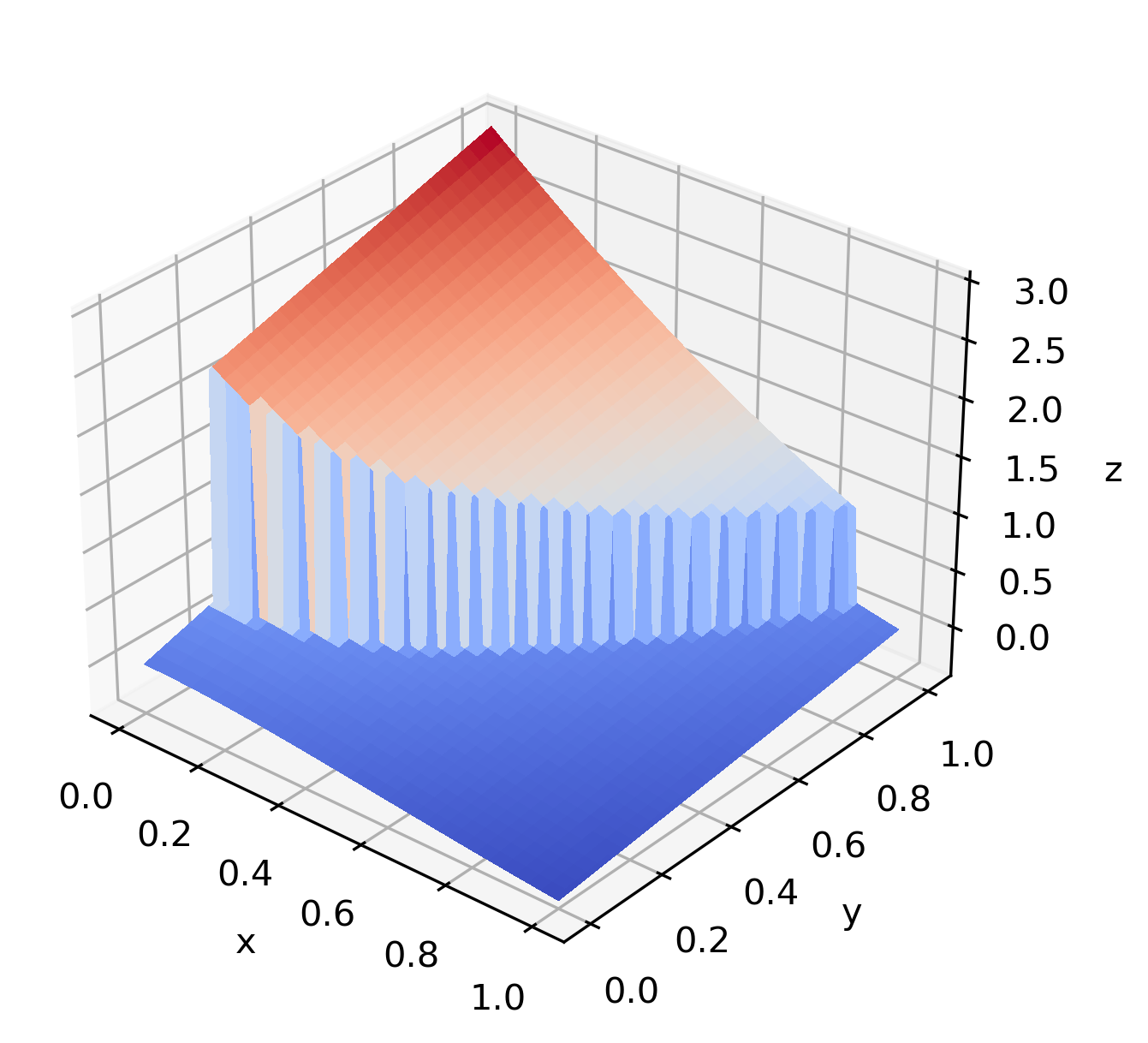}
\end{minipage}%
}%
\hspace{0.2in}
\subfigure[A 2--60--60--1 ReLU NN function approximation\label{comparison5}]{
\begin{minipage}[t]{0.4\linewidth}
\centering
\includegraphics[width=2.1in]{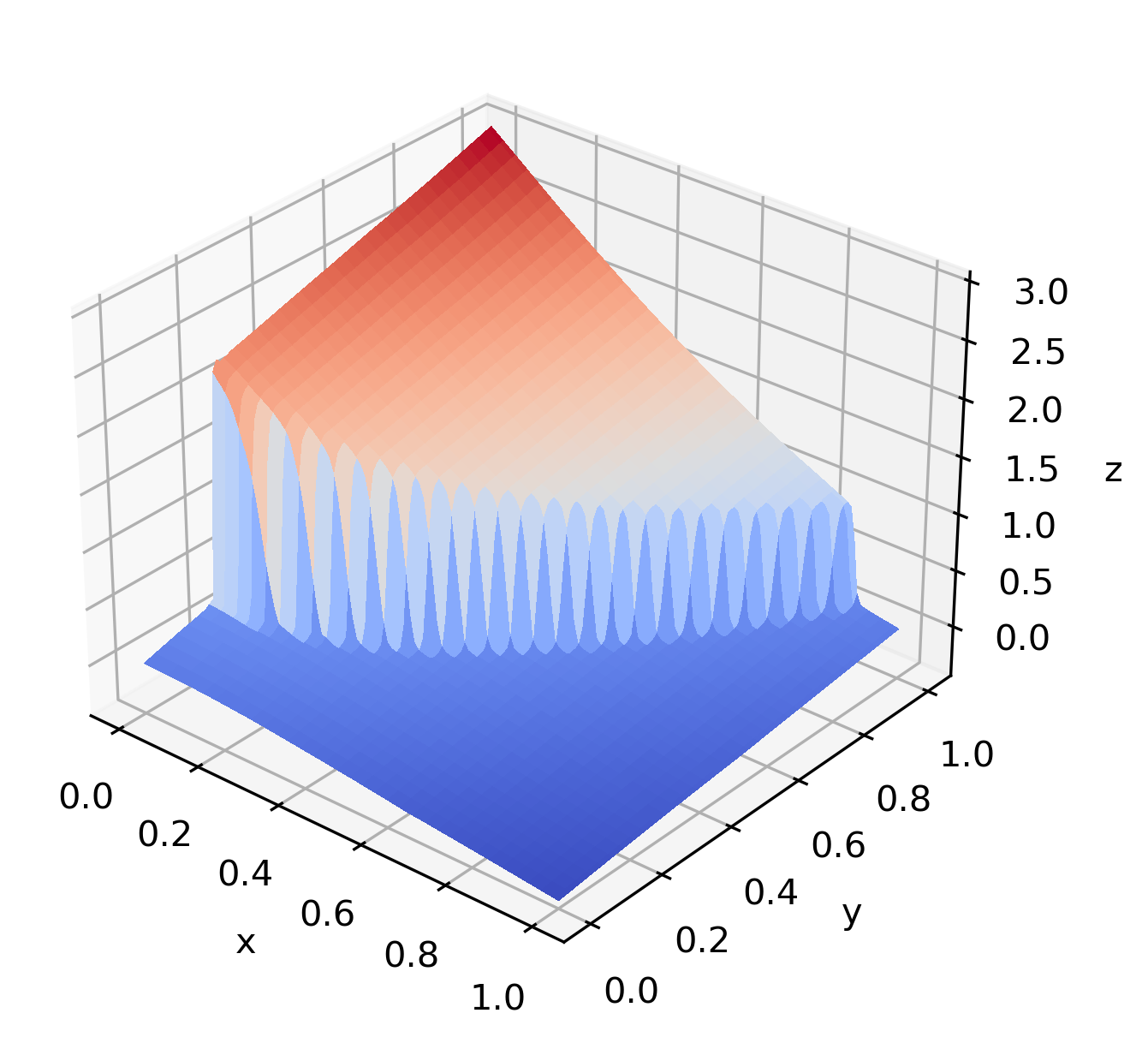}
\end{minipage}%
}%
\\
\subfigure[The breaking hyper-planes of the approximation in Figure \ref{comparison5}\label{breaking5}]{
\begin{minipage}[t]{0.4\linewidth}
\centering
\includegraphics[width=2.1in]{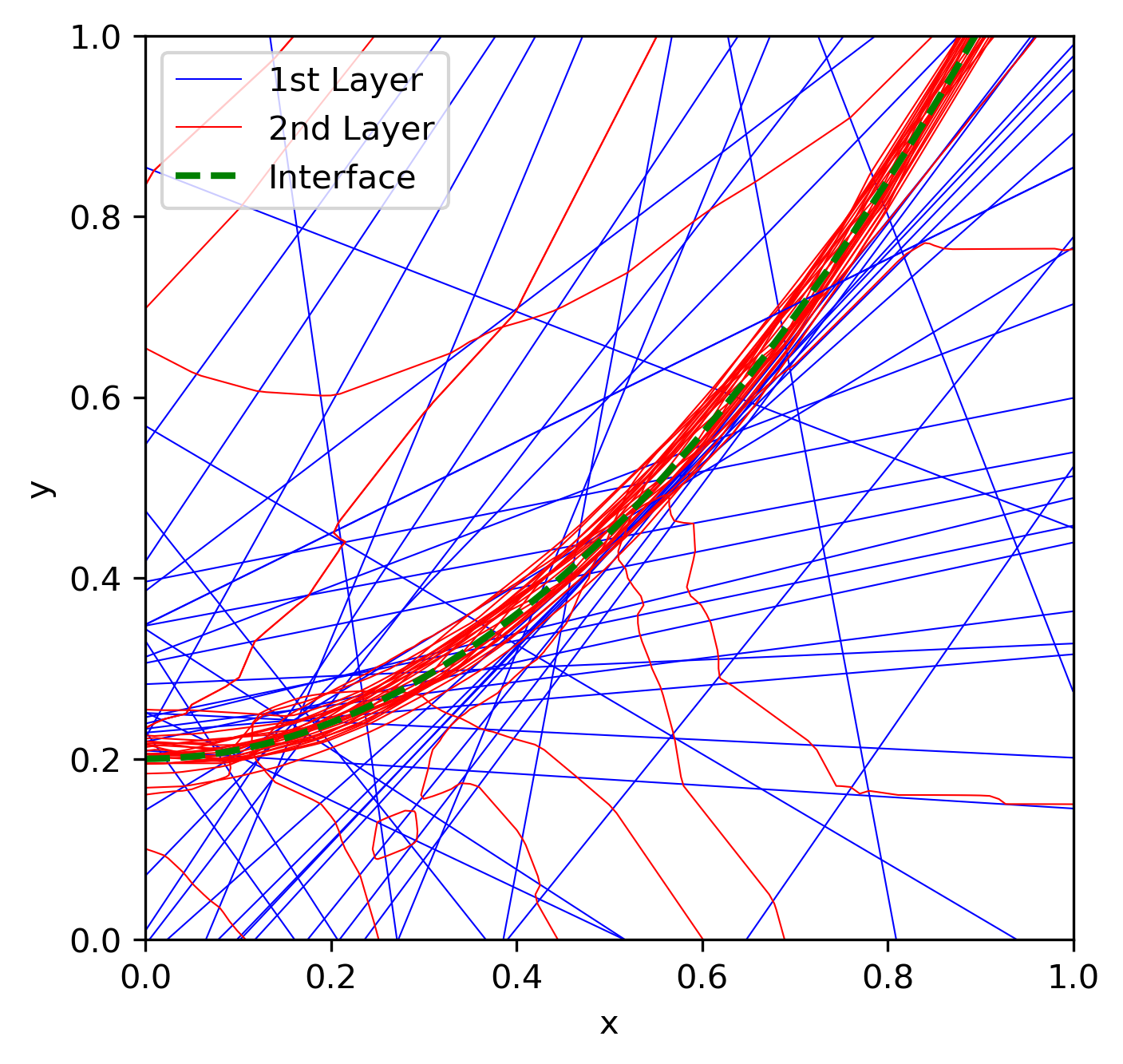}
\end{minipage}%
}%
\hspace{0.1in}
\subfigure[The loss curve ]{
\begin{minipage}[t]{0.5\linewidth}
\centering
\includegraphics[width=2.5in]{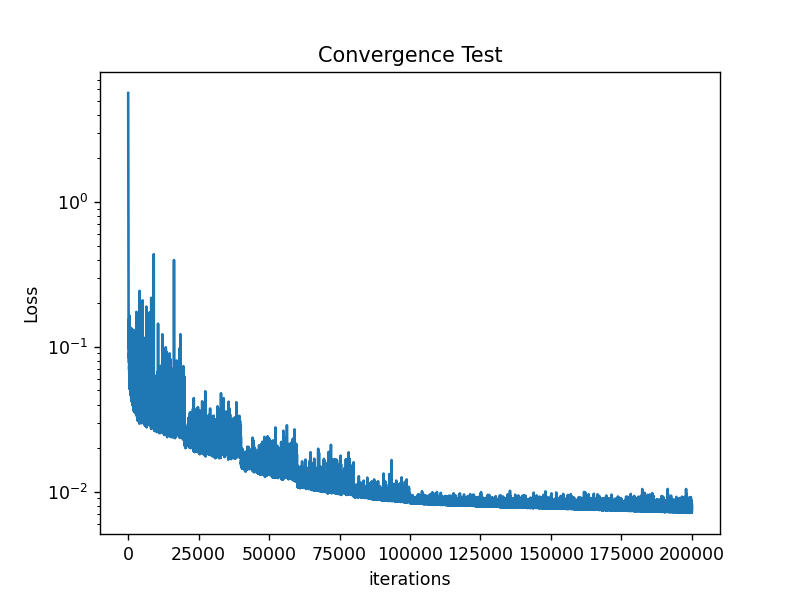}
\end{minipage}%
}%
\caption{Approximation results for the linear advection-reaction problem in Sec. 6.1.}
\end{figure}

The LSNN method was implemented using a 2--60--60--1 ReLU NN model and a uniform integration grid of size $h=0.01$ (see \cite{Cai2024nonconst} for experiment details). The directional derivative $v_{\bm\beta}$ was approximated by the backward finite difference quotient (\ref{finite_diff}) with $\rho=h/10$. We report the numerical results after 200,000 Adam iterations in \cref{test5 figure,test5 table}. As shown in \cref{vertical5,comparison_exact5,comparison5}, the LSNN method is capable of approximating the discontinuous solution with the curved interface with non-constant jump accurately without any oscillation or overshooting. In \cref{breaking5}, the graph of the physical mesh created by the trained ReLU NN function shows that the optimization process tends to distribute the breaking polylines in the second layer along the interface (see \cref{interface5}) presented in the problem, allowing the discontinuous solution to be accurately approximated using a piecewise linear function. \cref{test5 table} shows the relative numerical errors measured in different norms. With 3841 parameters, the ReLU NN can accurately approximate the solution with reasonable accuracy. 

\begin{table}[htbp]\label{test5 table}
\caption{Relative errors of the linear advection-reaction problem.}
\centering
\begin{tabular}{|l|l|l|l|l|}
\hline
Network structure  &$\frac{\|u-{u}_{_{\tiny N,\cT}}\|_0}{\|u\|_0}$ &$\frac{\vertiii{u-{u}_{_{\footnotesize N,\!\cT}}}_{\bm\beta}}{\vertiii{u}_{\bm\beta}}$ & $\frac{\mathcal{L}^{1/2}\left({u}_{_{\tiny N,\cT}};{\bf g}\right)}{\mathcal{L}^{1/2}\left({u}_{_{\tiny N,\cT}};{\bf 0}\right)}$ & Parameters \\[3mm] \hline
2--60--60--1  & 0.071953 & 0.115680 & 0.035981   & 3841\\ \hline
\end{tabular}
\end{table}

\begin{table}[htbp]
\centering
\caption{Relative $L^2$ errors of LSNN for the Riemann problem with $f(u)=\frac{1}{4}u^4$}
\vspace{5pt}
\begin{tabular}{|c|c|c|c|c|}
\hline
\multicolumn{2}{|c}{\multirow{2}{*}{Time block}} & \multicolumn{3}{|c|}{Number of sub-intervals} \\ \cline{3-5} 
\multicolumn{2}{|c|}{}& $\hat{m}=\hat{n}=2$  & $\hat{m}=\hat{n}=4$  & $\hat{m}=\hat{n}=6$  \\ \hline
\multirow{2}{*}{ $\Omega_{0,1}$} & Trapezoidal rule  &0.067712 &0.010446 &0.004543  \\ \cline{2-5} 
 &Mid-point rule&0.096238 &0.007917 &0.003381  \\ \hline
\multirow{2}{*}{ $\Omega_{1,2}$} & Trapezoidal rule  &0.108611 &0.008275 &0.009613  \\ \cline{2-5} 
 & Mid-point rule & 0.159651 &0.007169 &0.005028  \\ \hline
\end{tabular}
\label{shcl1_convergence}
\end{table}

% \begin{table}[htbp]
% \centering
% \caption{Relative $L^2$ errors of the problem with $f(u)=\frac{1}{4}u^4$ using the composite mid-point rule}
% \vspace{5pt}
% \begin{tabular}{|c|c|c|c|}
% \hline
% \multirow{2}{*}{Time block} & \multicolumn{3}{c|}{Number of sub-intervals} \\ \cline{2-4} 
%                             & $\hat{m}=\hat{n}=2$  & $\hat{m}=\hat{n}=4$  & $\hat{m}=\hat{n}=6$  \\ \bottomrule
% \multicolumn{1}{|l|}{ $\Omega_{0,1}$} &0.096238 &0.007917 &0.003381  \\ \hline
% \multicolumn{1}{|l|}{ $\Omega_{1,2}$} &0.159651 &0.007169 &0.005028  \\ \hline
% \end{tabular}
% \label{midpoint_convergence}
% \end{table}

\begin{figure}[ht]
\centering
    \subfigure[Traces at $t=0.2$ %of exact and numerical solutions $u_{1,_\cT}$ 
    (trapezoidal) %on the plane 
    ]{ 
    \includegraphics[width=2.2in]{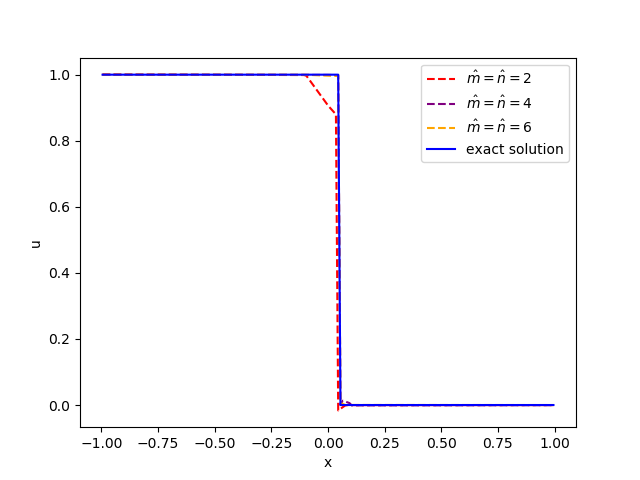}}
    \hspace{0.2in}
  \subfigure[Zoom-in plot near the discontinuous interface of sub-figure (a)]{
    \includegraphics[width=2.2in]{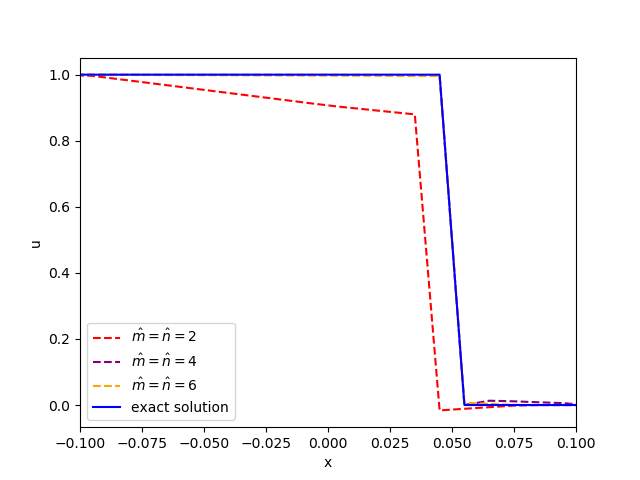}} \\
     % \hspace{0.1in}
  \subfigure[Traces at $t=0.4$ (trapezoidal) % rule on the plane 
  ]{
    \includegraphics[width=2.2 in]{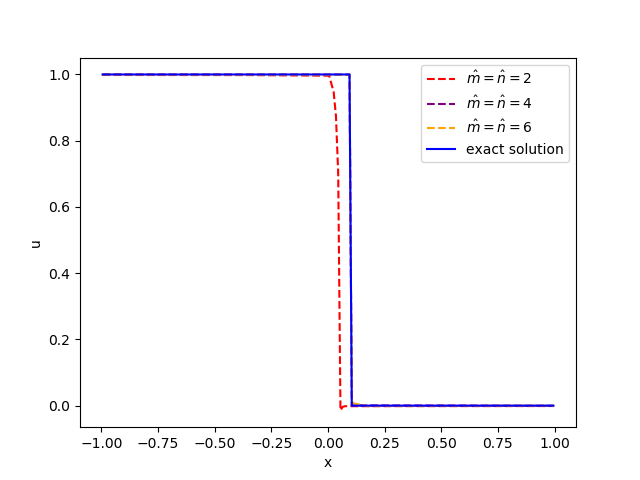}}
    \hspace{0.2in}
    \subfigure[Traces at $t=0.2$ (mid-point) ]{ 
    \includegraphics[width=2.2in]{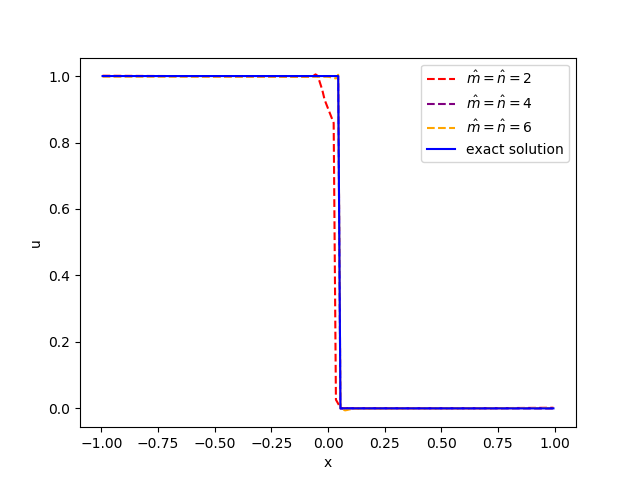}}\\
     % \hspace{0.1in}
  \subfigure[Zoom-in plot near the discontinuous interface of sub-figure (d)]{
    \includegraphics[width=2.2in]{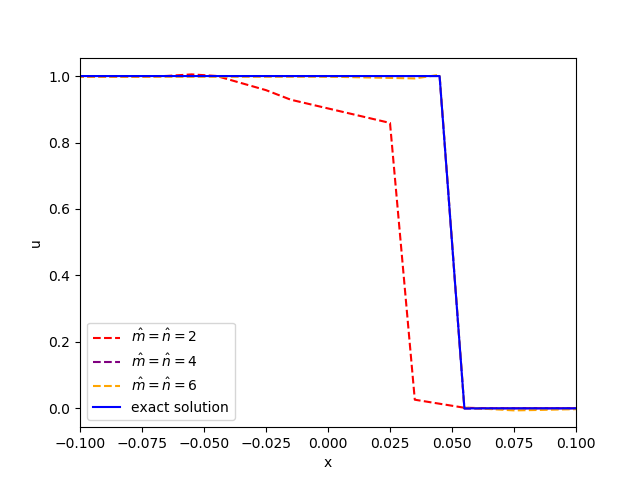}}
     \hspace{0.2in}
  \subfigure[Traces at $t=0.4$ (mid-point)]{
    \includegraphics[width=2.2in]{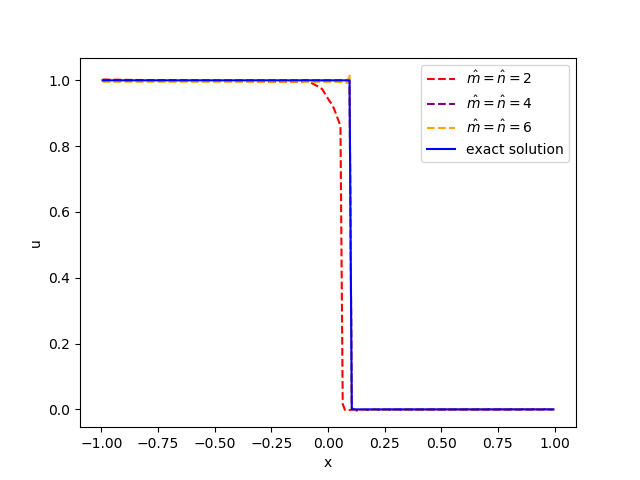}}
  \caption{Numerical results of the problem with $f(u)=\frac{1}{4}u^4$ using the composite trapezoidal and mid-point rules} 
  \label{convex_flux_fig}
\end{figure}

\subsection{A 1D Riemann problem with a spatial flux $f(u)=\frac14 u^4$}\label{sec6.2}
% The goals of this set of numerical experiments are twofold. First, we compare the performance of the LSNN method using 
% the composite trapezoidal/mid-point rule in (\ref{integration}). Second, we investigate the 

The second numerical example is a Riemann problem with a convex flux ${\bff}(u) = (f(u), u)=(\frac14 u^4, u)$ and an initial condition with a unit jump $u_{_L}=1 >0=u_{_R}$ at $(0,0)$ \cite{Cai2023nonlinear}. 
The computational domain is chosen as $\Omega = (-1,1)\times (0,0.4)$ and is subdivided into two blocks, $\Omega_{0,1} = (-1,1)\times (0,0.2)$ and $\Omega_{1,2} = (-1,1)\times (0.2,0.4)$, allowing efficient computation. The numerical integration is performed using a uniform grid of size $h_x=h_t=0.01$ and for the discrete divergence operator \divt \ref{div-dis2}, we tested two quadrature methods for calculating the line integral $\cQ_{\partial K}(\cdot)$: the composite trapezoidal rule and the midpoint rule. Furthermore, we also investigated the impact of the number of sub-intervals, along each boundary edge of $\partial K$, on the precision of the LSNN method.

A 2--10--10--1 ReLU NN model was used as an approximate function, and the Adam optimizer trains its associated parameters in $50,000$ iterations, the resulting relative $L^2$ errors are reported in Tables~\ref{shcl1_convergence}.And the traces of the exact and numerical solutions in  $t=0.2$ and $t=0.4$ are plotted in Fig.~\ref{convex_flux_fig}.

From Table~\ref{shcl1_convergence}, it is evident that the accuracy of the LSNN method depends on the number of sub-intervals, with $\hat{m}$ and $\hat{n}$ denoting the number of sub-intervals along the spatial and temporal directions, respectively. In general, larger values of $\hat{m}$ and $\hat{n}$ lead to higher accuracy in the LSNN approximation. However, increasing $\hat{m}$ and $\hat{n}$ also raises the computational cost, particularly for evaluating the line integral $\cQ_{\partial K}(\cdot)$, which becomes more expensive with finer partitions. Moreover, the accuracy achieved using the composite trapezoidal and midpoint rules in the LSNN method is comparable; both are capable of accurately simulating this Riemann problem and capturing the correct shock propagation speed.

% From Tables~\ref{shcl1_convergence}, it is expected that the accuracy of the LSNN method depends on the number of sub-intervals ($\hat{m}$ and $\hat{n}$ are the corresponding number of sub-intervals along the spatial and temporal directions, respectively); that is, the larger the $\hat{m}$ and $\hat{n}$, the more accurate the LSNN method is. Moreover, the accuracy using the composite trapezoidal and mid-point rules in the LSNN method is comparable, both are capable of simulating this Riemann problem with accurate shock propagating speed. 

\subsection{A 2D inviscid Burgers equation}%\label{sec_burgers}
The last numerical test considers a two-dimensional inviscid Burgers equation,  where the spatial flux vector field is
$\tilde{\bff}(u)=\frac12 (u^2, u^2)$. Given a piecewise constant initial data,
\begin{equation*}\label{2dburgers}
    u_0(x,y)=\left\{\begin{array}{r l}
    -0.2, & \mbox{if } \; x < 0.5\; \mbox{ and }\;y >0.5,\\
    -1.0, & \mbox{if }\; x > 0.5\; \mbox{ and }\; y >0.5,\\
    0.5, & \mbox{if } \; x < 0.5\; \mbox{ and } \; y < 0.5,\\
    0.8, & \mbox{if } \; x > 0.5\; \mbox{ and }\; y < 0.5,
    \end{array}\right.
\end{equation*}
the exact solution to this problem is as follows \cite{GUERMOND2014},
\begin{comment}
\begin{equation*}\label{2dburgers}
    u(x,y,t)=\left\{\begin{array}{r l}
    \begin{array}{cl}
   \left\{\begin{array}{ll} -0.2, & y>0.5 + 3t/20, \\ 0.5, & \mbox{otherwise} \end{array} \end{array}\right. & \mbox{and } \; x < 0.5-3t/5\; %\mbox{ and }\;
    %\left\{\begin{array}{ll}
    %y>0.5 + 3t/20,   \\ \mbox{otherwise}, \end{array}\right.   
    \\[4mm]
    \begin{array}{cl}
    -1 \\ 0.5 \end{array}  & \mbox{if } \; 0.5 - 3t/5 < x < 0.5-t/4\; \mbox{ and }\;
    \left\{\begin{array}{ll}
    y>-8x/7+15/14-15t/28,   \\ \mbox{otherwise},
    \end{array}\right.   \\
    \begin{array}{cl}
    -1 \\ 0.5 \end{array}  & \mbox{if } \; 0.5 - t/4 < x < 0.5+t/2\; \mbox{ and }\;
    \left\{\begin{array}{ll}
    y>x/6+ 5/12-5t/24,   \\ \mbox{otherwise},
    \end{array}\right.   \\
    \begin{array}{cl}
    -1 \\ \frac{2x-1}{2t} \end{array}  & \mbox{if } \; 0.5 + t/2 < x < 0.5 + 4t/5\; \mbox{ and }\;
    \left\{\begin{array}{ll}
    y> x-\frac{5}{18t}(x+t-0.5)^2,   \\ \mbox{otherwise},
    \end{array}\right.   \\
    \begin{array}{cl}
    -1 \\ 0.8 \end{array}  & \mbox{if } \; x > 0.5+ 4t/5\; \mbox{ and }\;
    \left\{\begin{array}{ll}
    y>0.5-t/10,   \\ \mbox{otherwise}.
    \end{array}\right.   \\
    \end{array}\right.
\end{equation*}
\end{comment}
\begin{equation*}\label{2dburgers_sol}
    u(x,y,t)=\left\{\begin{array}{l l}
    %\begin{array}{cl}
   \left\{\begin{array}{ll} -0.2, & y>0.5 + 3t/20, \\[1mm] 0.5, & \mbox{otherwise} \end{array}\right. & \mbox{and } \; x < 0.5-3t/5, %\mbox{ and }\;
    %\left\{\begin{array}{ll}
    %y>0.5 + 3t/20,   \\ \mbox{otherwise}, \end{array}\right.   
    \\[4mm]
    \left\{\begin{array}{ll}
    -1, & y>-8x/7+15/14-15t/28,\\[1mm] 0.5, & \mbox{otherwise}  \end{array}\right.  
    & \mbox{and } \; 0.5 - 3t/5 < x < 0.5-t/4, %\mbox{ and }\;\left\{\begin{array}{ll} y>-8x/7+15/14-15t/28,   \\ \mbox{otherwise}, \end{array}\right.   
    \\[4mm]
    \left\{\begin{array}{ll}
    -1, & y>x/6+ 5/12-5t/24, \\[1mm] 0.5, & \mbox{otherwise} \end{array}\right.  & \mbox{and } \; 0.5 - t/4 < x < 0.5+t/2,
    %\left\{\begin{array}{ll}   y>x/6+ 5/12-5t/24,   \\ \mbox{otherwise},    \end{array}\right.   
    \\[4mm]
    \left\{\begin{array}{ll}
    -1, & y> x-\frac{5}{18t}(x+t-0.5)^2, \\[1mm] \frac{1}{2t} (2x-1), & \mbox{otherwise}\end{array}\right.  & \mbox{and } \; 0.5 + t/2 < x < 0.5 + 4t/5, %\mbox{ and }\;\left\{\begin{array}{ll}y> x-\frac{5}{18t}(x+t-0.5)^2,   \\ \mbox{otherwise},\end{array}\right.   
    \\[4mm]
    \left\{\begin{array}{cl}
    -1, & y>0.5-t/10, \\[1mm] 0.8, & \mbox{otherwise} \end{array}\right.  & \mbox{and } \; x > 0.5+ 4t/5. 
    %\left\{\begin{array}{ll}    y>0.5-t/10,   \\ \mbox{otherwise}.    \end{array}\right.   \\
    \end{array}\right.
\end{equation*}

\begin{figure}[htbp]
\centering
    \subfigure[$t=0.1$]{ 
    \includegraphics[width=2in]{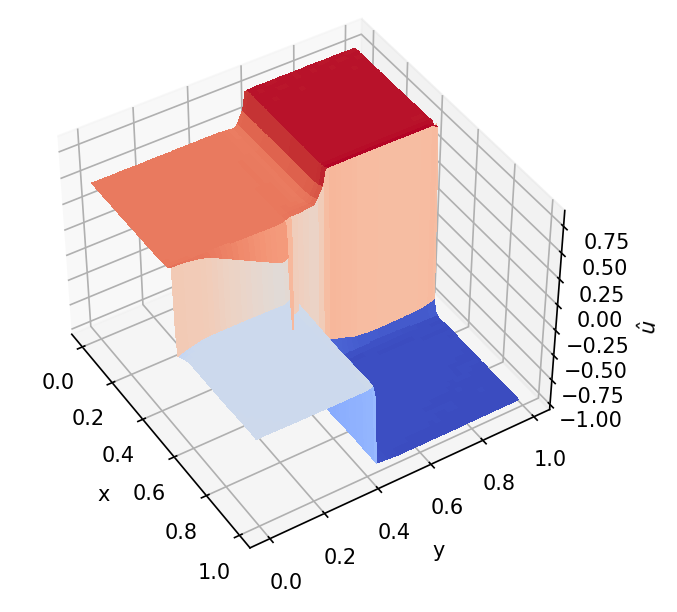}}
  \subfigure[$t=0.3$ ]{
    \includegraphics[width=1.8in]{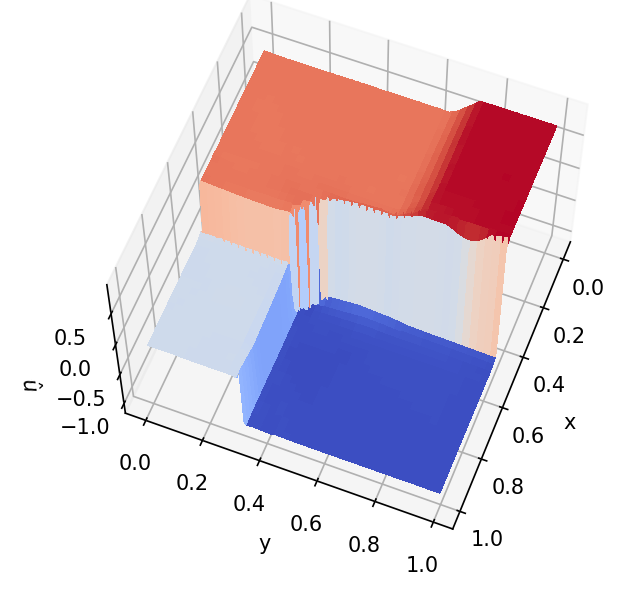}}
    \subfigure[$t=0.5$]{ 
    \includegraphics[width=1.8in]{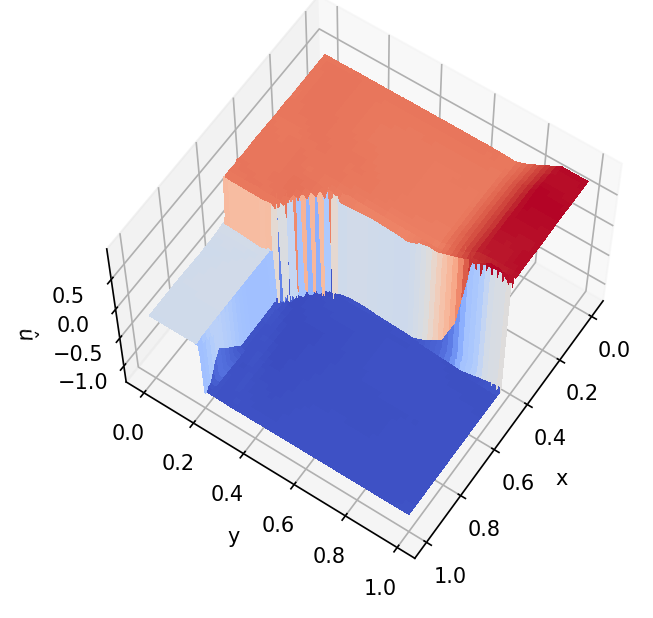}}
  \caption{Numerical results of $2D$ Burgers' equation.} 
  \label{fig_burger2d}
\end{figure}

Setting the computational domain $\Omega = (0,1)^2\times (0, 0.5)$, and the inflow boundary conditions prescribed using the exact solution, a 4-layer ReLU NN (3--48--48--48--1) was used as the model function. Again, the numerical integration was performed on uniform grids of size $h_x=h_y=h_t=0.01$, and the computation domain is decomposed into five time blocks of equal sizes, namely $\Omega_{0,1}, \Omega_{1,2}, \cdots, \Omega_{4,5}$. The three-dimensional discrete divergence operator \divt is computed using the mid-point quadrature rule with $\hat{m}=\hat{n}=\hat{k}=2$, where $\hat{m}$, $\hat{n}$ and $\hat{k}$ are the number of sub-intervals along the spatial $x$, spatial $y$ and the temporal direction. Table \ref{riemann_2d_table} reported the relative $L^2$errors of LSNN in each time block. Specifically, $30,000$ iterations of Adam optimization were performed for the first time block, and the rest blocks were trained with $20,000$ iterations. Fig.\ref{fig_burger2d} presents the numerical results at time $t= 0.1$, $0.3$, and $0.5$. This experiment shows that the LSNN method can be extended to two-dimensional problems and is capable of simulating the shock and rarefaction waves simultaneously.

As anticipated, numerical error accumulated when using a block space-time method mentioned in the previous paragraph that decomposes the time interval $[0,5]$ into $5$ blocks (see \cite{Cai2021nonlinear}). By $t= 0.5$, the relative error $L^2$ reached $21.3\%$ (see Table \ref{riemann_2d_table}) . This result raises an important question for future research: how to enhance the accuracy of the LSNN method for high-dimensional hyperbolic problems. Theoretical studies suggest that a three-layer ReLU NN is sufficient for such problems from a function approximation standpoint \cite{CCL2024}. However, developing an efficient and reliable iterative solver suitable for these high-dimensional, non-convex optimization problems remains a challenge. The discussion in Sec. \ref{s:solver} offers insights into leveraging the unique structure of NNs to guide the iterative process, though the problem remains unresolved.%maybe need to say something to address the accumulated error issue. 

\begin{table}[htbp]
\centering
\caption{Relative $L^2$ errors of LSNN for a 2D Burgers' equation}
\vspace{5pt}
\begin{tabular}{|l|l|l|}
\hline
Network structure &Block & $\frac{\|u^k-u^k_{_\cT}\|_0}{\|u^k\|_0}$ \\ \hline
 \multirow{3}{*}{3-48-48-48-1} & $\Omega_{0,1}$  & 0.093679 \\ \cline{2-3}
 & $\Omega_{1,2}$  & 0.121375 \\ \cline{2-3}
 & $\Omega_{2,3}$  & 0.163755\\ \cline{2-3}
 & $\Omega_{3,4}$  & 0.190460\\ \cline{2-3}
 & $\Omega_{4,5}$  & 0.213013\\ \hline
 
\end{tabular}
\centering
\label{riemann_2d_table}
\end{table}

% \section{Conclusion}

%\section{Code availability}

%The code used in the study will be publicly available through the GitHub.

\newpage
\bigskip
\bibliographystyle{ieee}
\bibliography{Reference}

\end{document}

%% file: ex_shared.tex
\usepackage{lipsum}
\usepackage{amsfonts}
\usepackage{graphicx}
\usepackage{epstopdf}
\usepackage{algorithmic}
\ifpdf
  \DeclareGraphicsExtensions{.eps,.pdf,.png,.jpg}
\else
  \DeclareGraphicsExtensions{.eps}
\fi

\newsiamremark{hypothesis}{Hypothesis}
\crefname{hypothesis}{Hypothesis}{Hypotheses}
\newsiamthm{claim}{Claim}

\headers{Least-Squares Neural Network Method}{M. Liu and Z. Cai}

\title{Least-Squares Neural Network (LSNN) Method \\ [2mm]  for Scalar Hyperbolic Partial Differential Equations} %\thanks{This work was supported in part by the National Science Foundation under grant DMS-2110571.}}

\author{Min Liu\thanks{School of Mechanical Engineering, Purdue University, 585 Purdue Mall,
West Lafayette, IN 47907-2088 (\email{liu66@purdue.edu}). }\and Zhiqiang Cai\thanks{Department of Mathematics, Purdue University, 150 N. University Street, West Lafayette, IN 47907-2067 
  (\email{caiz@purdue.edu}).}
}

\usepackage{amsopn}